\documentclass[10pt, a4paper, english, twocolumn, final]{extarticle}

\usepackage{adjustbox} % Nice alternative to minipage.
\usepackage{afterpage} % To give the title page its own geometry.
\usepackage[boxed,noline,linesnumbered,resetcount,algosection,noend]{algorithm2e}
\usepackage{algpseudocode}

\usepackage{amsmath} % Nice maths symbols.
\usepackage{amssymb} % Nice variable symbols.
\let\amsDiamond=\Diamond
\usepackage{amsthm}
\usepackage{array} % Allow for custom column widths in tables.
\usepackage{ltablex} % For long tables spanning multiple pages. % Must be before ARYDSHLN package!
\keepXColumns % Keeps the X column
\usepackage{arydshln} % Dashed lines using \hdashline \cdashline
\usepackage[english]{babel} % (I only use this to give blindmath math symbols.)
\usepackage{bbm} % Gives Blackboard fonts.
\usepackage{blindtext} % Generates dummy maths. cf. lipsum.
\usepackage{calc} % Calculates widths of words. 

\usepackage{chngcntr} % Changing counters, e.g. with footnotes.
\usepackage{cuted} % For single column modes in twocolumn mode. 
\usepackage{cprotect} % For verbatim like environments in todonotes. 
\usepackage{datenumber}
\usepackage{datetime}
\usepackage[scale=3, color={[gray]{0.9}}]{draftwatermark} % Gives a draft overlay. Use options [nostamp] or [final].
\usepackage[spaced=-1,mleftright]{diffcoeff}
\usepackage{emptypage} % Empty pages have no headers and footers.
\usepackage{paralist} % Must be loaded before enumitem.
\usepackage{enumitem} % Nice listing options in itemize and enumerate. 
\usepackage{etoolbox} % For defining conditionals. 
\usepackage{fancyhdr} % Nice headers.
\usepackage{float} % Nice figure placement.
\usepackage[T1]{fontenc} % Nice range of text characters and accents.
\usepackage[bottom,multiple]{footmisc} % Nice footnote formatting.
\usepackage{graphicx} % Include figures.
\usepackage[notquote]{hanging} % For indenting later lines in a paragraph. USE the 'noquote' option else the `'` is overwritten, and breaks in maths mode!

\usepackage{ifoddpage} % Checks for odd or even page.
\usepackage[geometry]{ifsym} % Useful symbols.
\usepackage{imakeidx} % Makes the index.
\usepackage{indentfirst} % Indents the first paragraph.
\usepackage{letltxmacro} % For defining a nice SQRT symbol.
\usepackage{lineno} % For line numbers. 
\usepackage{lipsum} % Useful for adding jargon.
\usepackage{listings} % The listings package for code.
\usepackage{marginnote} % For nice margin notes.
\usepackage{mathtools} % Gives the colon equals symbol.
\usepackage[framed,numbered,useliterate]{mcode} % Inports Listings package ideal for MATLAB.
\usepackage[framemethod=tikz]{mdframed} % Gives nices boxed and sidesrules.
\usepackage{mleftright}
\usepackage{multirow} % Nice table cells spanning many rows.
\usepackage{multicol} % If I want to use multiple columns.
\usepackage[numbers, sort&compress]{natbib} % Nice references.
\usepackage{bm} % Bold math symbols. %% Needs to be loaded after newpxtext/math. %%
\usepackage{nicefrac} % Gives nice fractions for superscripts.
\usepackage{nomencl} % Gives a symbol nomenclature. 
\usepackage[super]{nth} % Gives nice ordinal superscripts, eg 1st, 2nd, etc.
\usepackage{parskip} % Gives nicer indenting.
\usepackage{ragged2e} % For nice allignment.
\usepackage{realboxes}
\usepackage{romannum} % Nice typing for roman numerals.
\usepackage{rotating} % For sideways figures.
\usepackage[scr,scaled=1.1]{rsfso} % Gives Script fonts which are not so slanted. 
\usepackage{setspace} % Ideal for increasing line spacing. E.g.  \doublespacing
\usepackage{siunitx} % Nice formating of units.
\usepackage{soul}
\usepackage{subcaption} % Side by side figures.
\usepackage[textsize=small]{todonotes} % A nice TODO list. [disable] to supress.
\usepackage{tikz} % Nice diagrams.
\usepackage{titling} % Access title variables. 
\usepackage{titlesec} % Nice section title colouring options.
\usepackage[nottoc]{tocbibind} % Gives nices Table of Contents
\usepackage[hyphens,spaces,obeyspaces]{url}
\usepackage{verbatimbox} % For verbatim inside todonotes. 
\usepackage{wasysym}
\usepackage{xcolor} % This is useful for making greyed table cells, nice for headers. Known preamble placement issues.
\usepackage{xifthen}% Provides \isempty test.
\usepackage{xparse} % Gives \NewDocumentEnvironment which has nice optional argument handling.
\usepackage{xpatch}
\usepackage{xspace} % Gives nice spacing for commands.
\usepackage{xurl} % Nicer leanbreaks for URLs. 
\usepackage[colorlinks=true,linkcolor=black,urlcolor=black,citecolor=black,anchorcolor=black]{hyperref} % Colour links.
\usepackage[noabbrev]{cleveref} % Gives smart referencing. %% After Hyperref
\usepackage[margin=10pt,font=small,labelfont=bf,labelsep=endash,figurewithin=section,tablewithin=section]{caption} % Caption figures and tables nicely. %% After cleveref.
\usepackage[margin=20mm]{geometry} % Use nice margins. Does give a small change in the default page margins. 
\usepackage{graphicx}

\makeindex

\makenomenclature
\makeatletter
\g@addto@macro\bfseries{\boldmath}
\makeatother

\preto\subequations{\ifhmode\unskip\fi}

\SetAlCapSkip{1em} % Increase algorithm spacing between box and caption.
\SetAlCapNameFnt{\small}
\SetAlCapFnt{\small}
\SetAlFnt{\small}
\SetAlgoCaptionSeparator{ \textendash}
\makeatletter
\patchcmd{\@algocf@start}% <cmd>
  {-1.5em}% <search>
  {-1em}% <replace>
{}{}% <success><failure>
\makeatother
\newcommand{\algrule}[1][0.2pt]{\par\vskip.2\baselineskip\hrule height #1\par\vskip.2\baselineskip}
\SetKw{KwGoTo}{go to}
\SetKwIF{If}{ElseIf}{Else}{if}{}{else if}{else}{end if}%
\SetKwFor{While}{while}{}{endw}%
\SetKwFor{ForEach}{for each}{}{endfch}
\SetKwRepeat{Do}{do}{while}
\SetInd{1em}{0em}

\makeatletter
\renewrobustcmd{\cref}{\@osmcref{cref}}
\renewrobustcmd{\Cref}{\@osmcref{Cref}}
\def\@osmcref#1#2{%
    \begingroup
    \ifcsundef{r@#2}
    {}
    {\expandafter\expandafter\expandafter\expandafter\expandafter
        \expandafter\expandafter\def
        \expandafter\expandafter\expandafter\expandafter\expandafter
        \expandafter\expandafter\@osmcref@name
        \expandafter\expandafter\expandafter\expandafter\expandafter
        \expandafter\expandafter{%
            \expandafter\expandafter\expandafter
            \@thirdoffive\csname r@#2\endcsname}}%
    \ifcsundef{r@#2@cref}
    {}
    {\cref@gettype{#2}{\@osmcref@type}}%
    \ifboolexpr{not test {\ifdefvoid{\@osmcref@name}}
        and (test {\ifdefstring{\@osmcref@type}{theorem}}
        or test {\ifdefstring{\@osmcref@type}{lemma}}
        or test {\ifdefstring{\@osmcref@type}{corollary}})}
    {\nameref{#2} (\@cref{#1}{#2})}
    {\@cref{#1}{#2}}%
    \endgroup
}
\makeatother

\newtheorem{theorem}{Theorem}[section]
\newtheorem{corollary}{Corollary}[theorem]

\theoremstyle{definition}
\newtheorem{remark}{Remark}[section]
\newtheorem{definition}{Definition}[section]
\newtheorem{assumption}{Assumption}[section]
\crefname{lemma}{lemma}{lemmas}
\Crefname{lemma}{Lemma}{Lemmas}
\crefname{theorem}{theorem}{theorems}
\Crefname{theorem}{Theorem}{Theorems}
\crefname{corollary}{corollary}{corollaries}
\Crefname{corollary}{Corollary}{Corollaries}
\Crefname{proposition}{proposition}{propositions}
\Crefname{proposition}{Proposition}{Propositions}
\crefname{remark}{remark}{remarks}
\Crefname{remark}{Remark}{Remarks}
\crefname{definition}{definition}{definitions}
\Crefname{definition}{Definition}{Definitions}
\crefname{assumption}{assumption}{assumptions}
\Crefname{assumption}{Assumption}{Assumptions}
\newcommand{\prooffont}{\upshape\bfseries}
\xpatchcmd{\proof}{\itshape}{\prooffont}{}{}

\makeatletter
\def\th@plain{%
    \thm@notefont{}% same as heading font
    \itshape % body font
}
\def\th@definition{%
    \thm@notefont{}% same as heading font
    \normalfont % body font
}
\makeatother

\graphicspath{{.}}

\crefname{listing}{code}{codes}
\Crefname{listing}{Code}{Codes}
\definecolor{codegray}{rgb}{0.9,0.9,0.9}
\newfloat{lstfloat}{htbp}{lop} % environment for placing lisings in to make them float. 

\addto\captionsenglish{} % When using babel
\addto\captionsenglish{} % When using babel

\newcolumntype{L}[1]{>{\raggedright\let\newline\\\arraybackslash\hspace{0pt}}m{#1}}
\newcolumntype{C}[1]{>{\centering\let\newline\\\arraybackslash\hspace{0pt}}m{#1}}
\newcolumntype{R}[1]{>{\raggedleft\let\newline\\\arraybackslash\hspace{0pt}}m{#1}}

\makeatletter
\let\oldr@@t\r@@t
\def\r@@t#1#2{%
    \setbox0=\hbox{\(\oldr@@t#1{#2\,}\)}\dimen0=\ht0
    \advance\dimen0-0.2\ht0
    \setbox2=\hbox{\vrule height\ht0 depth -\dimen0}%
    {\box0\lower0.4pt\box2}}
\LetLtxMacro{\oldsqrt}{\sqrt}
\renewcommand*{\sqrt}[2][\ ]{\oldsqrt[#1]{#2}}
\makeatother

\DeclareMathOperator*{\argmin}{argmin}

\numberwithin{equation}{section}
\newcommand*\tageq{\refstepcounter{equation}\tag{\theequation}}
\crefname{equation}{}{}
\creflabelformat{equation}{(#2#1#3)}

\newcommand{\firstrowspacing}{\rule{0pt}{2.6ex}}
\setlist{listparindent=\parindent,parsep=1ex} 

\def\sator{Sator Arepo tenet opera rotas.\xspace}

\newcount\loopcounter
\def\dummysentences#1{%
    \loopcounter = #1
    \loop
    \sator\ %
    \advance\loopcounter by -1
    \ifnum\loopcounter > 0
    \repeat%
}

\blindmathtrue

\NewCommandCopy\oldtodo\todo
\outer\def\todo[#1]{\icprotect{\oldtodo[#1]}}

\makeatletter
\xpretocmd\lstinline{\Colorbox{codegray}\bgroup\appto\lst@DeInit{\egroup}}{}{}
\makeatother

\lstdefinestyle{longline}
{numbers=none,
frame=none, 
xleftmargin=0.2em, 
xrightmargin=0em, 
framexleftmargin=0.2em,
}

\newcommand{\mylistlabelfont}[1]{\normalfont\textit{#1:}}
\newlist{longdescription}{description}{1}
\setlist[longdescription]{%
  style=unboxed,
  font=\mylistlabelfont
}

\DeclareMathOperator{\trace}{trace}
\DeclareMathOperator{\diag}{diag}
\xpatchcmd{\mathdisplay}{$$}{$$\textstyle}{}{}

\title{Convergence proofs and strong error bounds for forward-backward stochastic differential equations using neural network simulations}
\author{Dr Oliver Sheridan-Methven\thanks{\href{mailto:oliver.sheridan-methven@hotmail.co.uk}{\texttt{oliver.sheridan-methven@hotmail.co.uk}}.\newline \hspace*{0.8\footnotemargin} The code used to generate the results and figures in this report is hosted at \href{https://github.com/oliversheridanmethven/mlmc_with_nn_dissertation}{\url{https://github.com/oliversheridanmethven/mlmc\textunderscore with\textunderscore nn}}.}}
\date{}
\begin{document}
%\linenumbers
\pagenumbering{arabic}

\setlength{\droptitle}{-3em}   % To raise or lower the maketitle
\maketitle
%\input{contents}

%\doublespacing
\newpage
\begin{longdescription}
\item[Keywords] Forward-backward stochastic differential equations, neural networks, stochastic simulation, multilevel Monte Carlo, Feynman-Kac theorem, and strong error bounds. 
\item[MSC codes] 65C05, 65C30, 60H35, 93E03, 65M15, 68U20, 68W40, 93E99, and 60G99. 
\end{longdescription}
\begin{abstract}
We introduce forward-backward stochastic differential equations, highlighting the connection between solutions of these and solutions of partial differential equations, related by the Feynman-Kac theorem. We review the technique of approximating solutions to high dimensional partial differential equations using neural networks, and similarly approximating solutions of stochastic differential equations using multilevel Monte Carlo. Connecting the multilevel Monte Carlo method with the neural network framework using the setup established by \citet{e2017deep_learning,e2018solving,e2021algorithms} and \citet{raissi2018forwardbackward}, we provide novel numerical analyses to produce strong error bounds for the specific framework of \citet{raissi2018forwardbackward}. Our results bound the overall strong error in terms of the maximum of the discretisation error and the neural network's approximation error. Our analyses are necessary for applications of multilevel Monte Carlo, for which we propose suitable frameworks to exploit the variance structures of the multilevel estimators we elucidate. Also, focusing on the loss function advocated by \citet{raissi2018forwardbackward}, we expose the limitations of this, highlighting and quantifying its bias and variance. Lastly, we propose various avenues of further research  which we anticipate should offer significant insight and speed improvements.
\end{abstract}
\section{Introduction}
\label{sec:introduction}

Continuously evolving systems driven by random underlying processes are found frequently in both academic and real world settings. Such systems are described by stochastic differential equations, and have found modelling applications in: physics, finance, statistics, material science, biology, chemistry, seismology, weather forecasting, robotics, automated vehicles, systems control, insect population outbreaks, neurology, epidemiology, criminology, urban planning, and numerous other settings. For some illustrative examples and a wide catalogue of references, we recommend the reader to: \citet[\S\,11--14]{klebaner2012introduction}, \citet[\S\,5.3, 18.3, and 20]{higham2021introduction}, and \citet[\S\,6--7]{kloeden1999numerical}. Associated to and extending upon these are forward-backward stochastic differential equations, whereby numerous stochastic processes are coupled together and must satisfy a combination of both initial conditions and terminal conditions, all the while being causally consistent (i.e.\ solutions are correctly adapted and use only present or previous information, but no future information). 

Solving and simulating such stochastic systems, either exactly or approximately, is an area of both great practical importance and academic interest. In recent decades and years computational capabilities have increased, techniques have improved, applications have become more ambitious, data have become abundant, and dimensionalities have become huge. All of these factors combine to make the state-of-the-art tasks very formidable. One particular difficulty of those just mentioned has been largely surmounted, which is that of high dimensionality. In recent years there has been an explosion in the use of various machine learning techniques to tackle problems with incredibly high dimensionalities, most prominently through the use of neural networks and similar tools. For works reviewing and surveying the current landscape we recommend the reader to \citet[\S\,10]{james2023intro}, \citet{abiodun2018state}, \citet{aloysius2017review}, and \citet{tripathi2010high}. 

The solutions of partial differential equations are directly related to the solutions of forward-backward stochastic differential equations, and \textit{vice versa}, through variants of the Feynman-Kac theorem. Neural networks have demonstrated their groundbreaking ability to accurately and efficiently approximate solutions to high dimensional partial differential equations. Consequently, this gives an entry point to using neural networks to simulate approximate solutions to high dimensional forward-backward stochastic differential equations, whose approximation is otherwise usually not readily accessible nor straightforward in general. 
Neural networks thus enable the simulation of forward-backward stochastic differential equations, and the simulation is subsequently done by Monte Carlo methods, with multilevel Monte Carlo methods being especially desirable because of their better rates of convergence \citep{giles2008multilevel}. These applications---where neural networks are combined with multilevel Monte Carlo---are what we concern ourselves with. A good and concise review of this problem setup is provided by \citet{e2018solving}.  

The core theme of this body of research investigates the interplay of training neural networks alongside generating sample path approximations. The primary concern is trying to ensure we are using an appropriately trained neural network at all stages in our calculations. At any given discretisation level used in the multilevel Monte Carlo setup, the neural network used should be: 
\begin{inparaenum}[1)]
  \item sufficiently well trained, but \label{item:well_trained}
  \item not trained more than is necessary. \label{item:not_overly_trained}
\end{inparaenum}
Unfortunately we find that \eqref{item:well_trained} is commonly assumed (or hoped for), and that \eqref{item:not_overly_trained} is not considered (or falsely written off when categorised as offline training). Our research focuses on developing a framework to ensure verifiably that both these criteria are satisfied, and providing corrective measures if they are not. The most related works in this direction are by \citet{e2017deep_learning,e2018solving,e2021algorithms} and \citet{raissi2018forwardbackward}, with peripheral works by \citet{guler2019robust} and \citet{naarayan2024thesis}. 

\subsection{Contributions of this report}

Our research has generated the following contributions, which we list in the order of their significance:
\begin{longdescription}
\item [Strong error bounds] 
We have shored up the empirical findings of \citet{raissi2018forwardbackward} with our own numerical analysis of their framework. We have shown many of the usual convergence orders found in regular stochastic systems carry forward to systems of coupled forward-backward stochastic differential equations utilising neural network approximations. In doing so we have both provided the supporting mathematical analysis, and highlighted what further assumptions and restrictions arise to guarantee such results. This mathematical underpinning which we provide to compliment and complete previously only empirical results is crucial for ensuring the soundness of the empirical results, and establishing them mathematically in the wider context of similar research in the field. 

\item [Loss function quantification] 
We critically inspect the loss function widely used by \citet{raissi2018forwardbackward}, \citet{e2017deep_learning,e2018solving,e2021algorithms}, and others, and make explicit a quantification of its bias and variance. This appears to have been broadly overlooked previously, whereas our treatment announces the existing limitations, and the consequences this has for the neural network's ability to approximate solutions to high dimensional partial differential equations. Consequently, our treatment suggests a simple modification to the loss function that should reduce its variance by a factor of \( \Delta t^{1/2} \) if Hessians are not prohibitively expensive to compute. 

\item [Multilevel Monte Carlo frameworks] 
We have presented a multilevel Monte Carlo framework to utilise differing neural network sophistications at differing temporal discretisations, and showcased the variance structure of the multilevel correction. 

\item [Highlight further avenues of research] We have given expositions of various ideas for further research, focusing on developing improved loss functions through a combination of their analytic properties, through conventional antithetic variance reduction techniques, and through differing interpolation points. 
\end{longdescription}

The significant and novel contributions herein are the error bounds, loss function analysis, and multilevel Monte Carlo framework.

\subsection{Structure of this report}

The remains of this report are structured as follows:
\begin{longdescription}
\item[\Cref{sec:mathematical_preliminaries}] 
Overviews all the mathematical preliminaries necessary for a technical description of the problem setup and the existing frameworks we will subsequently use. 

\item[\Cref{sec:numerical_analysis}]
Provides our numerical analysis and supporting experiment results for strong error bounds for path approximations generated using neural networks with the framework of \citet{raissi2018forwardbackward}, and associated results for appropriate multilevel Monte Carlo frameworks.

\item[\Cref{sec:future_research}] Elucidates avenues for further research, focusing on the loss function, antithetic techniques, and interpolation points. 

\item[\Cref{sec:conclusions}] Discusses the conclusions of this report. 

\end{longdescription}

\section{Mathematical preliminaries}
\label{sec:mathematical_preliminaries}

In this \namecref{sec:mathematical_preliminaries} we will briefly introduce the reader to forward-backward stochastic differential equations and overview the relation of solutions to these stochastic systems to the solutions of partial differential equations. With this relation established, we will review the applicability of neural networks for finding approximate solutions to the partial differential equations, completing the problem setup. This will bring us to the core contribution of this work, which will be investigations of multilevel Monte Carlo methodologies to provide further speed improvements to our framework, without compromising accuracy. 

\subsection{Forward-backward stochastic differential equations}
\label{sec:fbsdes}

We introduce the system of forward-backward stochastic differential equations, whereby we have the forward process
\begin{subequations}
\label{eqt:fbsde}
\begin{align}
\dl{X_t} & = a(t,X_t,Y_t,Z_t) \dl{t} + b(t, X_t, Y_t) \dl{W_t} \label{eqt:fbsde_fsde}\\
\shortintertext{with \( X_{t} = X_0 \) at \(t = 0\), and the backward process}
\dl{Y_t} & = \phi(t,X_t,Y_t,Z_t) \dl{t} + Z_t^\top \dl{W_t}  \label{eqt:fbsde_bsde}
\end{align}
\end{subequations}
with \( Y_T = \xi(X_T) \), where \( Y_t \) is c\`{a}dl\`{a}g and adapted and \( Z_t \) is predictable. In lieu of a conventional name for this \(Z_t\) process, we hereby refer to it as the \textit{hidden} process. Notice that we do not have a stochastic differential equation to describe the evolution of \(Z_t\). Rather, finding \(Z_t\) is considered part of the problem to be solved (e.g.\ analogous to determining free boundaries when solving partial differential equations).

Forward-backward stochastic differential equations such as these frequently find applications in e.g.\ mathematical finance and stochastic control \citep{peng1993bsde_and_applications}, but also in physics, ecology, neuroscience, etc.
\citet[p.\,219]{gobet2016monte} further highlights the connection to stochastic control, and for a more thorough resource linking stochastic control and forward-backward stochastic differential equations we recommend \citet[\S\,2]{nusken2021solving}.

Here \( X_t \in \mathbb{R}^d \), and thus we have a \( d \)-dimensional stochastic process. Clarifying the dimensions, domains, and ranges for the various terms, we have:
\begin{itemize}
\item \( X_t, W_t, Z_t \in \mathbb{R}^d \).
\item \( Y_t \in \mathbb{R} \).
\item \( a \colon [0, T]\times \mathbb{R}^d\times \mathbb{R} \times \mathbb{R}^{d\times d} \to  \mathbb{R}^d\).
\item \( b \colon [0, T]\times \mathbb{R}^d\times \mathbb{R} \to  \mathbb{R}^{d\times d}\).
\item \( \xi \colon \mathbb{R}^d \to \mathbb{R} \).
\item \( \phi \colon [0, T]\times \mathbb{R}^d\times \mathbb{R} \times \mathbb{R}^{d\times d} \to  \mathbb{R}\).
\end{itemize}

To make any headway with \cref{eqt:fbsde}, we require a host of standard assumptions \citep[\S\,3.2, 4.5, and 10.2]{kloeden1999numerical}, which will enable us to posit the existence and uniqueness of solutions and appproximations thereof.

\begin{assumption}
\label{asmp:measurability}
The functions \( a \)  and \( b \)  are jointly (Lebesgue) \( \mathcal{L}^2 \)-measurable.
\end{assumption}

\begin{assumption}
\label{asmp:lipschitz_continuity}
The functions \( a \)  and \( b \) are Lipschitz continuous in their spatial dimensions.
\end{assumption}

\begin{assumption}
\label{asmp:linear_growth_bound}
The functions \( a \)  and \( b \) have linear spatial growth.
\end{assumption}

\begin{assumption}
\label{asmp:measurable_initial_condition}
\( X_{0} \) is measurable with \( \mathbb{E}(\lVert X_{0} \rVert^2_2) < \infty\).
\end{assumption}

\Cref{asmp:measurability,asmp:lipschitz_continuity} are sufficient to ensure uniqueness of any solutions to \cref{eqt:fbsde_fsde} for regular stochastic differential equations \citep[lemma~4.5.2]{kloeden1999numerical}, and similarly \cref{asmp:measurability,asmp:lipschitz_continuity,asmp:linear_growth_bound,asmp:measurable_initial_condition} ensure there exists a solution \citep[theorem~4.5.3]{kloeden1999numerical}. To ensure the convergence of numerical approximations \( \widehat{X} \) of the solutions using uniform time steps \( \Delta t \) (such as e.g.\ the Euler-Maruyama scheme \citep{kloeden1999numerical}), we require additional assumptions \citep[\S\,10.2]{kloeden1999numerical}.

\begin{assumption}
\label{asmp:temporal_holder_continuity}
The functions \( a \)  and \( b \) are \(\tfrac{1}{2}\)-H\"{o}lder continuous with linear spatial growth.
\end{assumption}

\begin{assumption}
	\label{asmp:convergent_initial_value}
	There exists a constant \( K > 0 \) such that \( \mathbb{E}(\lVert X_{0} - \widehat{X}_{0}\rVert^2_2) \leq K\Delta t \).
\end{assumption}

To assert the existence of solutions to \cref{eqt:fbsde_bsde} for regular backward stochastic differential equations, we require further standard assumptions \citep[\S\,19.1]{cohen2015stochastic}.

\begin{assumption}
	\label{asmp:bsde_standard_data}
	The data \( (\xi, \phi) \) are \textit{standard}, whereby \( \mathbb{E}(\xi^2) < \infty \) and \( \mathbb{E}(\int_{0}^T \phi^2(t, \cdot, \cdot, \cdot) \dl{t}) < \infty \).
\end{assumption}

\begin{assumption}
	\label{asmp:bsde_standard_lipschitz_data}
	The data \( (\xi, \phi) \) are \textit{standard Lipschitz}, whereby they satisfy \cref{asmp:bsde_standard_data} and \( \phi \) is spatially Lipschitz.
\end{assumption}

If \( (\xi, \phi) \) satisfy \cref{asmp:bsde_standard_lipschitz_data} then \cref{eqt:fbsde_bsde} with data \( (\xi, \phi) \) for regular backward stochastic differential equations admits a unique solution \( (Y_t, Z_t) \) \citep[theorem~19.1.7]{cohen2015stochastic}, and for any stopping time \( s \in [0, T] \) with \(s \geq t\) almost surely the pair \( (Y_t, Z_t) \) is also the unique solution to \cref{eqt:fbsde_bsde} with data \( (Y_s, \phi) \) \citep[lemma~19.1.8]{cohen2015stochastic}.

\begin{assumption}
\label{asmp:invertible_and_bounded_diffusion}
\citep[p.\,221]{gobet2016monte}
The diffusion process \( b \) is everywhere invertible and the inverse is uniformly bounded.
\end{assumption}

Solutions of the fully coupled forward-backward stochastic differential equations \cref{eqt:fbsde} are connected to the solutions of partial differential equations through \cref{thm:feynman_kac_semi_linear_multidimensional}  \citep{peng2016bsde}. 

\begin{remark}
	To clarify our vector calculus notation and avoid any confusion or ambiguity we mirror the popular notation used by e.g.\ \citet{nocedal2006numerical}. Namely,  
	\( \nabla \) represents the gradient operator \( \nabla_i \coloneq \diffp{}{x_i} \) and \( \nabla^2 \) denotes the Hessian operator \( \nabla^2_{ij} \coloneq \diffp{}{x_i,x_j} \).\footnote{The Hessian operator is not to be confused with the equally common Laplacian operator \( \nabla^2 \coloneq \sum_{i = 1}^d \diffp[2]{}{x_i} \), which occurs very frequent in literature concerning partial differential equations \citep[p.\,352]{riley2010mathematical}.}
\end{remark}

\begin{theorem}[the multi-dimensional semi-linear Feynman-Kac theorem]
\label{thm:feynman_kac_semi_linear_multidimensional}
\citep[theorem~7.2.1]{gobet2016monte}
\citep{peng1991probabilistic}
\citep{pardoux1992backward}
\citep{wu2014probabilistic}
\citep[theorem~19.5.1]{cohen2015stochastic}
Let \( L \) be the differential operator
\begin{multline}
\label{eqt:feynman_kac_spatial_pde_operator}
L u(t,x)  \coloneq   \diffp{u}{t}(t, x) + \sum_{i = 1}^{d} \alpha_i \nabla_i u(t, x) \\
+ \frac{1}{2} \sum^{d}_{i,j=1} (\beta \beta^\top)_{ij} \nabla^2_{ij} u (t, x),
\end{multline}
where 
{\allowdisplaybreaks
\begin{align}
\alpha & \coloneq a(t, x, u(t, x), \nabla u (t, x)) \\
\shortintertext{and}
\beta & \coloneq b(t, x, u(t, x)),
\end{align}
}%
and where \( a \) and \( b \) are the drift and diffusion functions appearing in \cref{eqt:fbsde_fsde}.
Let \( u \colon [0, T] \times \mathbb{R}^d \to \mathbb{R} \) be the solution to the semi-linear partial differential equation 
\begin{equation}
\label{eqt:feynman_kac_semi_linear_pde}
L u(t, x) = \phi(t, x, u(t, x), \nabla u(t, x)) 
\end{equation}
with terminal boundary condition \( u(T, x) = g(x) \) for some known function \( g \),
which for some positive constant \( K \) satisfies 
\begin{equation}
\label{eqt:pde_solution_lipschitz_bounds}
\sup_{(t,x) \in [0, T]\times \mathbb{R}^d} \frac{\lvert u \rvert^2 + \lVert \nabla u \rVert^2_2 }{1 + \lVert x\rVert^2_2} \leq K.
\end{equation}

If a solution \( u \) exists, and if \( X_t \) is the solution to \cref{eqt:fbsde_fsde} with initial condition \( X_0 \), then 
\begin{equation}
(Y_t, Z_t) \coloneq (u(t, X_t), (b(t, X_t))^\top \nabla u (t, X_t)) \label{eqt:data_dependence_on_pde_solution}
\end{equation}
is the unique solution to \cref{eqt:fbsde_bsde} with data \( (g(X_T), \phi) \).
\end{theorem}

\begin{remark}
It is important to note the direction of implication in \cref{thm:feynman_kac_semi_linear_multidimensional}, and that solutions of the partial differential equation give rise to solutions of the forward-backward stochastic differential equations. The reverse implication does not necessarily follow, as discussed by \citet[remark~4.4]{karatzas1998brownian}. For a comprehensive discussion of Feynman-Kac formulas, we recommend \citet[\S\,4.4]{karatzas1998brownian}.
\end{remark}

The numerical approximation of solutions of general forward-backward stochastic differential equations by means of stochastic simulation techniques has no obvious entry point, because of a combination of the backward nature arising from the presence of a terminal condition, and the additional difficulty of ensuring the processes \( Y_t \) and \( Z_t \) are appropriately adapted and predictable. This makes directly simulating solutions of \cref{eqt:fbsde} a difficult and usually intractable task in general. However, in the special case where the forward stochastic process is conveniently decoupled from the backward stochastic process, then it is possible to simulate the two using e.g.\ a forward and backward Euler-Maruyama scheme. Examples of performing such simulations and the associated error analysis are given by \citet{massing2023approximation}, and earlier proposals of discretisation methods are also explored by \citet{zhang2004numerical}, and the related Monte Carlo simulations by \citet{bouchard2004discrete}. \citet{bouchard2004discrete} discuss the problem of ensuring the backwards Euler-Maruyama scheme produces a correctly adapted process, and possible modifications that recover this property. 

A further point to comment on is the difficulty of working in multiple dimensions. In going from a one dimensional forward stochastic process to a multidimensional one, not only does this place an added conceptual and theoretic difficulty, but also a computational one. The onset of the added difficulty though is present even in a low dimensional setting. Thus far we have been almost exclusively discussing the Euler-Maruyama scheme as our preferred simulation method. However, there do exist numerical schemes with higher convergence orders \citep[\S\,10]{kloeden1999numerical}, such as e.g.\ the Milstein scheme. Considering the Milstein scheme, in the multidimensional setting, the simulation of the forward stochastic process requires sampling L\'{e}vy areas \citep[p.\,343--344]{glasserman2003monte}, and this is in general a very difficult problem. Methods for sampling these have been developed by \citet{gaines1994random} and \citet{alnafisah2018implementation} (which are efficient in 2-dimensions), or for trying to circumvent them altogether, such as by exploiting antithetic twin paths by \citet{giles2014antithetic}.

\subsection{High dimensional partial differential equations}

Solving well behaved partial differential equations in one spatial dimension is usually easy, and there are countless software packages to tackle this. However, as the dimensionality increases, the computational difficulty significantly increases as well. Methods for solving well behaved partial differential equations in two or three spatial dimensions (and sometimes four) are well explored and remain a highly active area of research. A driving factor for this is the abundance of real world physical systems that exist in relatively low dimensions. 
However, underlying most methods for approximating solutions of partial differential equations are discretisation schemes, splitting the continuous spatial domain into a discrete mesh of points (whether they form a regular mesh or an irregular one), producing an exponential growth in the computational complexity, known as the curse of dimensionality \citep{seydel2006tools}. While Monte Carlo convergence is typically quite slow, it does not suffer from the same pitfall \citep[p.\,2--3]{glasserman2003monte}.

It is common place to utilise the relation between partial differential equations and stochastic processes, and to find the solution of either problem by means of solving the other. However, this is not the only means, and for example \citet{bender2008time_discretization} use a Picard iteration method for coupled forward-backward stochastic differential equations, which does not exploit their relation to partial differential equations. This has been extended more recently by \citet{e2019multilevel_picard,e2021multilevel_picard} who highlight this as a means of proving the approximation capabilities of deep neural networks, which \citet{hutzenthaler2020overcoming_pdes} showcase as a means of overcoming the curse of dimensionality \citep{hutzenthaler2020overcoming_pdes,hutzenthaler2020overcoming_pricing}.

\subsection{The Black-Scholes-Barenblatt equation}

We present the Black-Scholes-Barenblatt equation\footnote{The Black-Scholes-Barenblatt equation seems to have been named by \citet[\S\,1, p.\,76]{avellaneda1995pricing}, who to the best of our knowledge provide the first reference to the equation, crediting \citeauthor{barenblatt1979similarity} by stating: ``\textit{The physicist G.~I.~\citet{barenblatt1979similarity} ([sic]~1978) introduced a diffusion equation with a similar non-linearity to model flow in porous media; hence our terminology}''. Note that \citet{barenblatt1979similarity}~(\citeyear{barenblatt1979similarity}) is not to be confused with \citet{barenblatt1996similarity}~(\citeyear{barenblatt1996similarity}), where both \citep{barenblatt1979similarity} and \citep{barenblatt1996similarity} are [confusingly] printed by the same author with the same titles.} \citep{avellaneda1995pricing} with dimensionality 100 as an example of a forward-backward stochastic differential equation. Other examples include the Hamilton-Jacobi-Bellman equation \citep[\S\,4.3]{yong1999stochastic}, and the Allen-Cahn equation \citep[\S\,3, (12)]{allen1979microscopic}. These examples are taken from (and popularised by) \citet{e2017deep_learning,e2018solving} and \citet{raissi2018forwardbackward}, where several other examples are elucidated by \citet[\S\,4.4--4.7]{e2017deep_learning}.

A simplified constant volatility form of the Black-Scholes-Barenblatt partial differential equation \citep[(7)]{avellaneda1995pricing} is
\begin{equation}
\label{eqt:black_scholes_barenblatt_pde}
\diffp{u}{t} = r(u - (\nabla u)^\top x) - \frac{1}{2}\trace(\sigma^2\diag(X_t^2)\nabla^2 u)
\end{equation}
where \(r\) and \(\sigma\) are strictly positive constants. With terminal condition \(u(T, x) = g(x)\) this has the explicit solution
\begin{equation}
\label{eqt:black_scholes_barenblatt_pde_exact_solution}
u(t,x) = \exp((r+\sigma^2)(T-t))g(x).
\end{equation}
This is the solution (through \cref{thm:feynman_kac_semi_linear_multidimensional}) to the partial differential equation underlying the forward-backward stochastic differential equations
{\allowdisplaybreaks
\begin{subequations}
\begin{align}
\dl{X_t} & = \sigma \diag(X_t) \dl{W_t} \\
\shortintertext{and}
\dl{Y_t} & = r(Y_t - Z_t^\top X_t) \dl{t} + \sigma Z_t^\top \diag(X_t) \dl{W_t},
\end{align}
\end{subequations}
}%
with \(Y_T = g(X_T)\). The forward equation is decoupled from the backward one, and is a driftless geometric Brownian motion, the exact solution of which is well known \citep[\S\,4.4, (4.6)]{kloeden1999numerical} to be 
\begin{equation}
X_t = X_0 \exp(-\tfrac{1}{2}\sigma^2t + \sigma W_t).
\end{equation}

It is because the various drift and diffusion processes are well behaved and that we have closed form solutions for \(X_t\) and \(u\) that our numeric examples will focus on this example. Additionally, as the diffusion is multiplicative, rather than additive, we expect the Euler-Maruyama scheme to show a strong convergence order of \(\tfrac{1}{2}\) \citep[\S\,9.3 and 10.2]{kloeden1999numerical}. If the diffusion process were additive (such as for e.g.\ arithmetic Brownian motion), then the Euler-Maruyama scheme can exhibit uncharacteristically fast rates of strong convergence \citep[p.\,341--342]{kloeden1999numerical}, achieving a strong convergence order of 1 \citep[\S\,10.3, exercise 10.3.4]{kloeden1999numerical}.

The general Black-Scholes-Barenblatt equation is discussed by \citet{avellaneda1995pricing}, and properties of solutions are discussed by \citet{vargiolu2001existence} and \citet{li2019properties}. The Hamilton-Jacobi-Bellman equation is introduced by \citet[\S\,6.5]{kloeden1999numerical} and a comprehensive treatment is provided by \citet{yong1999stochastic}. Exact solutions (closed or semi-closed) to the Allen-Cahn equation are in general not known, and require numerical approximation (e.g.\ the branching diffusion method \citep{henrylabordere2014numerical}). For a detailed discussion of numerical techniques to approximate solutions to the Allen-Cahn equation, we recommend \citet[\S\,6.2--6.3]{bartels2015numerical}. A recent weak error analysis of the Allen-Cahn equation is provided by \citet{breit2024weak_error}, and a strong error analysis by \citet{becker2023strong}.

\subsection{Neural networks}

In recent years, neural networks have proved through their applications to be an excellent tool for tackling high dimensional problems with large numbers of parameters, notably deep neural networks. So too have they found applications in approximating the solutions to high dimensional partial differential equations, and a good recent review is provided by \citet{germain2022approximation}. Similarly, \citet{e2018solving} solve high dimensional partial differential equations with deep neural networks, and more generally  \citet{e2021algorithms} discuss algorithms for tackling high dimensional partial differential equations through machine learning techniques (such as neural networks).\footnote{The related forward-backward stochastic differential equations considered by \citet{germain2022approximation} and \citet{e2018solving,e2021algorithms} are for systems where the forward process is decoupled from the backward one.} For readers unfamiliar with neural networks, we recommend  \citet[\S\,11]{hastie2009elements} and \citet{gurney2018introduction}.

The premiss of the method is to note through \cref{thm:feynman_kac_semi_linear_multidimensional} the relation between the solution of the partial differential equation to the system of forward-backward stochastic differential equations, and then discretise the stochastic process, and lastly minimise the discretisation error observed over a batch of sample paths. 

To flesh this out, we discretise the stochastic processes \cref{eqt:fbsde} over \( N + 1 \) uniformly spaced time points \( t_0 < t_1 < \cdots < t_N \) with \( t_n \coloneq n\Delta t \) for \( n \in \{0,1,\ldots,N\} \) where \( \Delta t \coloneq \tfrac{T}{N} \), and we look for approximate solutions \( (\widehat{X}_n,  \widehat{Y}_n, \widehat{Z}_n) \approx (X_{t_n}, Y_{t_n}, Z_{t_n})\). The na\"{i}ve method of discretisation (appropriate for high dimensions) is the Euler-Maruyama discretisation
\begin{subequations}\label{eqt:fbsde_em}
\begin{align}
\widehat{X}_{n+1} & \coloneq \begin{multlined}[t] \widehat{X}_n
 + a(t_n,\widehat{X}_n,\widehat{Y}_n,\widehat{Z}_n) \Delta t \\
 + b(t_n, \widehat{X}_n, \widehat{Y}_n) \Delta W_n
\end{multlined} \label{eqt:fbsde_fsde_em} \\
\shortintertext{and}
\widehat{Y}_{n+1} & \coloneq \widehat{Y}_n +  \phi(t_n,\widehat{X}_n,\widehat{Y}_n,\widehat{Z}_n) \Delta t + \widehat{Z}_n^\top \Delta W_n.  \label{eqt:fbsde_bsde_em}
\end{align}
\end{subequations}
Notice that this Euler-Maruyama discretisation does not provide any means of time stepping the \( \widehat{Z}_n \) approximation.

Now if we have an approximation \( \hat{u} \approx u \) to the solution of the underlying partial differential equation \cref{eqt:feynman_kac_semi_linear_pde}, where \(\hat{u}\) is parametrised by some \(\theta\) such that \(\hat{u}(s,x) \equiv \hat{u}(s,x;\theta)\) produces the corresponding approximations \( (\widehat{X}_n^{\theta},  \widehat{Y}_n^{\theta}, \widehat{Z}_n^{\theta}) \), then using
\cref{thm:feynman_kac_semi_linear_multidimensional}, from an initial condition \( \widehat{X}_0 \) we can generate from \(u\) the correctly adapted initial approximations
{\allowdisplaybreaks
\begin{subequations}\label{eqt:process_initial_values}
\begin{align}
\widehat{Y}_0 & \coloneq u(0, \widehat{X}_0) \label{eqt:backward_process_initial_value} \\
\shortintertext{and}
\widehat{Z}_0 & \coloneq (b(0, \widehat{X}_0))^\top \nabla u (0, \widehat{X}_0), 
\label{eqt:hidden_process_initial_value} \\
\shortintertext{and similarly from \(\hat{u}\) generate} 
\widehat{Y}_0^{\theta} & \coloneq \hat{u}(0, \widehat{X}_0;\theta) 
\label{eqt:backward_process_initial_value_approx}  \\
\shortintertext{and}
\widehat{Z}_0^{\theta} & \coloneq (b(0, \widehat{X}_0))^\top \nabla \hat{u} (0, \widehat{X}_0;\theta). \label{eqt:hidden_process_initial_value_approx}
\end{align}
\end{subequations}
}%
Now suppose after \( n \) steps we have using \(u\) generated the approximations \( (\widehat{X}_n,  \widehat{Y}_n, \widehat{Z}_n) \) and similarly using \(\hat{u}\) the approximations \( (\widehat{X}_n^{\theta},  \widehat{Y}_n^{\theta}, \widehat{Z}_n^{\theta}) \) with \(\widehat{X}_0^{\theta} \coloneq \widehat{X}_0\), then \cref{thm:feynman_kac_semi_linear_multidimensional} does not indicate how to generate \( \widehat{X}_{n+1} \) nor \( \widehat{X}_{n+1}^{\theta} \), so the only obvious possibility that remains is to use \cref{eqt:fbsde_fsde_em} for \( \widehat{X}_{n+1} \) and similarly for \( \widehat{X}_{n+1}^{\theta} \) use 
\begin{multline}
\widehat{X}_{n+1}^{\theta} \coloneqq \widehat{X}_n^{\theta} 
 + a(t_n,\widehat{X}_n^{\theta},\widehat{Y}_n^{\theta},\widehat{Z}_n^{\theta}) \Delta t \\
+ b(t_n, \widehat{X}_n^{\theta}, \widehat{Y}_n^{\theta}) \Delta W_n.
\label{eqt:fbsde_fsde_em_approx}
\end{multline}
Conversely, there is no Euler-Maruyama update scheme for \( \widehat{Z}_{n+1} \) nor \( \widehat{Z}_{n+1}^{\theta} \), so the most obvious remedy is to use 
{\allowdisplaybreaks
\begin{subequations}\label{eqt:hidden_process_forward_value_definitions}
\begin{align}
\widehat{Z}_{n+1} & \coloneq (b(t_{n+1}, \widehat{X}_{n+1}))^\top \nabla u (t_{n+1}, \widehat{X}_{n+1}) \label{eqt:hidden_process_forward_value_definition} \\
\shortintertext{and}
\widehat{Z}_{n+1}^{\theta} & \coloneq (b(t_{n+1}, \widehat{X}_{n+1}^{\theta}))^\top \nabla \hat{u} (t_{n+1}, \widehat{X}_{n+1}^{\theta};\theta) \label{eqt:hidden_process_forward_value_definition_approx} 
\end{align}
\end{subequations}
}%
from \cref{thm:feynman_kac_semi_linear_multidimensional}. However, for the \( \widehat{Y}_{n+1} \) estimate, we can either use \cref{eqt:fbsde_bsde_em} from the Euler-Maruyama scheme, or \cref{thm:feynman_kac_semi_linear_multidimensional}, and similarly for \( \widehat{Y}_{n+1}^{\theta} \).

The method proposed by \citet{raissi2018forwardbackward} and \citet{e2018solving} is to use \cref{thm:feynman_kac_semi_linear_multidimensional} and set 
{\allowdisplaybreaks
\begin{subequations}\label{eqt:backward_process_forward_value_definitions}
\begin{align}
\widehat{Y}_{n+1} & \coloneq u (t_{n+1}, \widehat{X}_{n+1}) \label{eqt:backward_process_forward_value_definition} \\
\shortintertext{and}
\widehat{Y}_{n+1}^{\theta} & \coloneq \hat{u} (t_{n+1}, \widehat{X}_{n+1}^{\theta};\theta)
\label{eqt:backward_process_forward_value_definition_approx}
\end{align} 
\end{subequations}
}%
Combining all these update schemes into a single method gives \cref{algo:euler_maruyama_with_nn}. 

\begin{algorithm}[htb]
	\DontPrintSemicolon
	\KwIn{An approximate solution \(\hat{u} \) and exact solution \(u\) to   \cref{eqt:feynman_kac_semi_linear_pde} and  approximate initial value \(\widehat{X}_0\) (which may be exact).}
	\KwOut{Approximate sample paths \((\widehat{X}, \widehat{Y}, \widehat{Z})\) and \((\widehat{X}^{\theta}, \widehat{Y}^{\theta}, \widehat{Z}^{\theta})\).}
	\algrule
	\SetKwProg{initialisation}{Initialisation}{}{}
	\initialisation{}{
        Set \(\widehat{X}_0^{\theta} \leftarrow \widehat{X}_0\).\;
		Set \(\widehat{Y}_0 \), \(\widehat{Z}_0 \), \(\widehat{Y}_0^{\theta} \), and \(\widehat{Z}_0^{\theta} \) using     \cref{eqt:process_initial_values}.\;
	}{}
	\SetKwProg{forwardupdates}{Forward updates}{}{}
	\forwardupdates{}{
		\For{\( n \in \{0, 1, \ldots, N-1\} \)}{
			Sample a Wiener increment \(\Delta W_n\).\;
			Set \( \widehat{X}_{n+1} \) using \cref{eqt:fbsde_fsde_em}.\;  
            Set \( \widehat{X}_{n+1}^{\theta} \) using \cref{eqt:fbsde_fsde_em_approx}.\;
			Set \(\widehat{Z}_{n+1} \) and \(\widehat{Z}_{n+1}^{\theta} \) using \cref{eqt:hidden_process_forward_value_definitions}.\;
			Set \(\widehat{Y}_{n+1} \) and \(\widehat{Y}_{n+1}^{\theta} \) using \cref{eqt:backward_process_forward_value_definitions}.\;
		}
		\Return \((\widehat{X}, \widehat{Y}, \widehat{Z})\) and \((\widehat{X}^{\theta}, \widehat{Y}^{\theta}, \widehat{Z}^{\theta})\).\;
	}{}
	\caption{The update scheme proposed by \citet{raissi2018forwardbackward}. If the exact solution \( u \) is unknown and not available, then all steps setting \((\widehat{X}, \widehat{Y}, \widehat{Z})\) values can be skipped, leaving only the approximation \((\widehat{X}^{\theta}, \widehat{Y}^{\theta}, \widehat{Z}^{\theta})\).}
	\label{algo:euler_maruyama_with_nn}
\end{algorithm}  

\begin{remark}
In \cref{algo:euler_maruyama_with_nn} the terminal value is still set using \( u \) and \( \hat{u} \), not the boundary value function \( g \), a fact reflected when constructing the loss function. 
\end{remark}

However, as \( \hat{u} \) is only an approximation, there is no guarantee that the \( \widehat{Y}^{\theta} \) process generated using \cref{eqt:backward_process_forward_value_definition_approx} satisfies the equivalent Euler-Maruyama discretisation equations \cref{eqt:fbsde_em}, for reasons discussed at length in \cref{sec:loss_function_of_raissi} (similarly neither would \( \widehat{Y} \) from \cref{eqt:backward_process_forward_value_definition}). The sentiment of \citet{raissi2018forwardbackward} (and also \citet{e2017deep_learning}) is that in some appropriate sense, minimising the errors in the Euler-Maruyama discretisation equations causes \( u \) to be learnt, and conversely learning \( u \) will minimise the errors in the Euler-Maruyama discretisation equations. This is implicitly posited heuristically by \citet{raissi2018forwardbackward}, whereas in \cref{sec:loss_function_of_raissi} we examine and quantify why this holds true, but also that this is limited by discretisation bias and variance.

The difference between the estimates produced by the two methods is hoped to act as a good proxy for the measure of the quality of the approximation \( \hat{u} \approx u \), and can consequently be used to construct a loss function. Minimising this loss function (over several sample path realisations of \( \widehat{X} \)) is how we learn a good approximation. 

The methods used by \citet{e2018solving} and \citet{han2020convergence} use a neural network to learn good approximations for the initial value of the backward process and the associated gradient, facilitating good approximations of the backward process' initial value.\footnote{\citet{han2020convergence} consider coupled forward-backward stochastic differential equations.} However, they limit themselves to just the initial value, whereas \citet{raissi2018forwardbackward} advocates learning a good approximation for the whole path, which has the added utility of then being able to better simulate entire path approximations for \( \widehat{Y} \). 

We will follow the more general approach advocated by \citet{raissi2018forwardbackward}, where our ambition is to be able to efficiently produce entire path approximations \( \widehat{Y} \). Consequently, we will use a single neural network to approximate \( \hat{u} \approx u \) over the entire temporal domain (unlike \citet{e2018solving} who use \( N-1 \) unique neural networks, one for each time point \( t_1, t_2, \ldots, t_{N-1} \)). We use \( \omega \) in the superscript to denote a specific Brownian motion realisation, where for example \( X^{(\omega)} \) a realisation of the forward process. Similarly for two different realisations \( \omega_1 \) and \( \omega_2 \) where \( \omega_1 \neq \omega_2 \), these will generate two different driving Brownian motions. For brevity, when writing \( \omega_m \)  for the \( m \)-th realisation we will omit the \( \omega \), writing just \( X^{(m)} \equiv X^{(\omega_m)}  \). Consequently, for a batch of \( M \) realisations, we mirror \citet[(6)]{raissi2018forwardbackward} and define the loss
\begin{align*}
\mathscr{L} & \coloneq  \sum_{m=1}^{M} \sum_{n=0}^{N-1} \lvert \widehat{Y}^{(m)}_{n+1}
\begin{aligned}[t]
& {}- \widehat{Y}^{(m)}_{n} - \phi^{(m)}_{n}\Delta t_n   \\
& {}- (\widehat{Z}^{(m)}_{n})^\top \Delta W^{(m)}_{n} \rvert^2
\end{aligned}
 \\
& \phantom{\coloneq} {} + \sum_{m=1}^{M} \lvert \widehat{Y}_N^{(m)} - g(\widehat{X}_N^{(m)}) \rvert^2,  \tageq \label{eqt:loss}
\end{align*}
where \( \phi^{(m)}_{n} \coloneq \phi(t_n, \widehat{X}^{(m)}_{n}, \widehat{Y}^{(m)}_{n}, \widehat{Z}^{(m)}_{n}) \). \citet{e2018solving,e2017deep_learning,e2021algorithms} approximate the initial condition and then use the Euler-Maruyama scheme \cref{eqt:fbsde_em} to time step to the final condition using just the final term as their loss function. \citet[p.\,5--6]{raissi2018forwardbackward} discusses at length the drawbacks from such an approach, and the more desirable facets of using \cref{eqt:loss} as the loss function. 

The implementation of the loss used by \citet{raissi2018forwardbackward} also includes an additional term, which is not discussed by \citet{raissi2018forwardbackward}, but is mentioned briefly by \citet{guler2019robust}, which is 
\begin{equation}
\sum_{m=1}^{M} \lVert \widehat{Z}_N^{(m)} - \nabla g(\widehat{X}_N^{(m)}) \rVert^2_2. \label{eqt:loss_extra_term}
\end{equation} 
This measures the terminal convergence of the hidden process, which is appropriate to include if \( g \) is differentiable. 

\begin{remark}
\label{remark:new_brownian_motion_samples_between_iterations}
In the implementation by \citet{raissi2018forwardbackward}, in each training iteration the loss function resamples a new batch of Brownian motions, and \citet{e2018solving} do exactly the same in their implementation. Whether this is intended, or an implementation oversight is unclear, as this detail is not discussed in \citep{raissi2018forwardbackward}. While this is not necessarily a bad idea, it also seems perfectly reasonable to keep the paths the same between each iteration. Whether this decision has a significant impact or not would be an interesting topic for further research, and would interleave well with the points of analysis discussed in \cref{sec:loss_function_of_raissi}. Certainly if we wish to train two neural networks and have them closely coupled (for multilevel Monte Carlo applications), then care needs to be taken that the underlying Brownian paths sampled are the same for each neural network at any iteration.
\end{remark}

\subsection{Multilevel Monte Carlo}
\label{sec:mlmc}

By employing neural networks, we have seen how we are able to generate approximate path simulations of the coupled backward stochastic process. However, the training of neural networks is expensive, as are evaluating them. Consequently, this raises concerns over the trade off between the quality of the approximations, and the  penality we pay for training and using high fidelity approximations. Multilevel Monte Carlo provides a framework for understanding and analysing this trade off, but moreover provides a route to constructing simulation frameworks which combine differing levels of speed and fidelity. The combination is such that mostly low cost crude calculations are combined telescopically with a few high cost high fidelity calculations to recover the high accuracy results. This gives the speed improvements of the low accuracy calculations, but maintains the accuracy of the high fidelity calculations. 

First developed for the two level setting by \citet{heinrich2001multilevel} for parametric integration, and extended to the general multilevel case for stochastic simulation by \citet{giles2008multilevel}, multilevel Monte Carlo provides a framework for combining estimators, capitalising on the high speed of crude estimators and the high accuracy of others. Suppose there are \( L+1 \) levels of estimators, where each level is indexed by some \( l \in {0,1,\ldots,L}\), where \( l = 0 \) is the crudest, and \( l = L \) the finest, and intermediate values of \( l \) form a spectrum between these two extremes. Similarly, suppose that the cruder estimators are very cheap to compute, whereas the finer estimators are expensive to compute. For a review of multilevel Monte Carlo, we recommend \citet{giles2015multilevel}.

We are purposefully using terms such as ``fine'' and ``crude'' rather nebulously. This is to reflect that whatever characteristic property differentiating finer and cruder levels of the estimators could be one (or many) of several differences. Readers more familiar with multilevel Monte Carlo will recognise one such means of differentiating the various levels is through the level of temporal discretisation, where the crude levels use a very coarse time step, and the fine levels use considerably finer ones (as is the setting considered by \citet{giles2008multilevel}). However, there are a variety of other characteristics that can be varied to achieve the same effect. For example, the crude levels could use random variables sampled from cheaper approximate distributions \citep{giles2022approximate,giles2023approximating}, or could also use lower floating point precisions \citep{giles2024rounding}. More general approaches (not using approximate random variables) also fall under the category of multi-index Monte Carlo, developed by \citet{haji2016multi}. 

Having the different levels, we wish to exactly estimate some property \( P \) of the underlying stochastic process of interest \( X_t \). As \( X \) typically must be approximated by some \( \widehat{X} \approx X \), we must compromise on using \( P(X_t) \) and use an estimator \( \widehat{P}(\widehat{X}_n) \approx P(X_{t_n}) \). Denoting \( \widehat{P}_l \) as the approximation of \( \widehat{P} \) using the level \( l \) framework,
the usual Monte Carlo and multilevel Monte Carlo decomposition can be written as
\begin{equation}
\label{eqt:mc_and_mlmc}
\mathbb{E}(P) \approx \mathbb{E}(\widehat{P}_L) = \mathbb{E}(P_0) + \sum_{l = 1}^{L} \mathbb{E}(\widehat{P}_{l} - \widehat{P}_{l-1}),
\end{equation}
where the first approximation is the usual Monte Carlo framework using the highest level \( L \) of fidelity, and the subsequent equality is the multilevel Monte Carlo framework using a collection of multiple levels. For notational convenience, we will adopt the convention that \( \widehat{P}_{-1} \coloneq 0 \), allowing us to write our multilevel Monte Carlo decomposition as 
\begin{equation}
\label{eqt:mlmc}
\mathbb{E}(\widehat{P}_L) = \sum_{l = 0}^{L} \mathbb{E}(\widehat{P}_{l} - \widehat{P}_{l-1}).
\end{equation}
If we have another means of further estimating our existing estimators, where \( \widetilde{P} \approx \widehat{P} \), then we can nest a multilevel Monte Carlo with another, giving rise to the nested multilevel Monte Carlo framework
\begin{equation}
\label{eqt:nested_mlmc}
\mathbb{E}(\widehat{P}_L) = \sum_{l = 0}^{L} 
\begin{multlined}[t]
\mathbb{E}(\widetilde{P}_{l} - \widetilde{P}_{l-1}) \\ - \mathbb{E}(\widehat{P}_{l} - \widehat{P}_{l-1} - \widetilde{P}_{l} + \widetilde{P}_{l-1}).
\end{multlined}
\end{equation}

\subsubsection{Coupling the levels}

For the multilevel Monte Carlo framework to be beneficial, we require the following conditions are met: 
\begin{inparaenum}[1)]
\item there are substantial time savings to be had when comparing the cost of sampling from the \( (l-1) \)-th level compared to the \( l \)-th, and \label{item:time_savings}
\item  the variance of the multilevel correction \( \widehat{P}_{l} - \widehat{P}_{l-1} \) be much less than the coarser estimator \( \widehat{P}_{l-1} \). \label{item:variance_reduction}
\end{inparaenum}
Of these, \eqref{item:time_savings} establishes the potential savings that can be achieved, and \eqref{item:variance_reduction} establishes the efficiency of recovering these savings. To ensure the lower variance of the multilevel estimator, the two levels must be coupled together somehow, and the core method of achieving this is having both levels use the same underlying Brownian motion paths when computing the multilevel correction. 

To see how this is done consider the simple stochastic differential equation
\begin{equation}
\dl{X_t} = a(t, X_t)\dl{t} + b(t, X_t) \dl{W_t},
\end{equation} 
for which we construct Euler-Maruyama approximations \(\widehat{X}\). In the temporal discretisation used by \citet{heinrich2001multilevel} and \citet{giles2008multilevel}, we construct two estimators, one using a fine temporal discretisation, and the other using a coarse one, denoted by \(\widehat{X}^{\mathrm{f}}\) and \(\widehat{X}^{\mathrm{c}}\) respectively. The coarse discretisation uses \(N\) equally sized time steps of size \(\Delta t \coloneq \tfrac{T}{N}\) such that the approximation is indexed by \(n \in \{0,1,\ldots,N\}\). For simplicity, we assume the finer discretisation is evaluated similarly, but using instead twice as many time steps,\footnote{\citet[\S\,4.1]{giles2008multilevel} discusses optimal choices for the sampling ratios between the coarse and finer levels, finding a factor of 7 to be the closest to optimal, although factors of 2 or 4, which are more computationally convenient, are still similarly competitive.} each of size \(\tfrac{\Delta t}{2}\), but is now indexed by \(n \in \{0, \tfrac{1}{2}, 1, 1+\tfrac{1}{2}, \ldots, N-\tfrac{1}{2}, N\}\). Consequently, the Euler-Maruyama approximation schemes for the two levels are
{\allowdisplaybreaks
\begin{subequations}\label{eqt:euler_maruyama_multilevel_discretisations}
\begin{align}
\widehat{X}_{n+1}^{\mathrm{c}} & \coloneq 
\begin{multlined}[t]
\widehat{X}_{n}^{\mathrm{c}} + a(n\Delta t, \widehat{X}_{n}^{\mathrm{c}}) \Delta t  \\ + b(n\Delta t, \widehat{X}_{n}^{\mathrm{c}}) \Delta W_n^{\mathrm{c}}
\end{multlined}
\label{eqt:euler_maruyama_coarse}\\
\shortintertext{and}
\widehat{X}_{n+1/2}^{\mathrm{f}} & \coloneq
\begin{multlined}[t]
\widehat{X}_{n}^{\mathrm{f}} + a(n\Delta t, \widehat{X}_{n}^{\mathrm{f}}) \tfrac{\Delta t}{2} \\ + b(n\Delta t, \widehat{X}_{n}^{\mathrm{f}}) \Delta W_n^{\mathrm{f}}
\end{multlined}
\label{eqt:euler_maruyama_fine}
\end{align}
\end{subequations}
}%
where \(\Delta W^{\mathrm{f}}_n \coloneq \sqrt{\Delta t / 2} Z_n\) with \(Z_n\) being independently and identically distributed standard Gaussian random variables, and the coarse paths being coupled to the fine path by enforcing 
\begin{equation}
\Delta W^{\mathrm{c}}_n \coloneq \Delta W^{\mathrm{f}}_n + \Delta W^{\mathrm{f}}_{n + 1/2}.
\end{equation}

\subsection{Training and inference}

Thus far, we have identified two distinct categories of work. The first is constructing the approximate solution of the underlying partial differential equation, and the second is using this in conjunction with the forward Euler-Maruyama scheme to generate sample paths. These correspond to training and inference stages respectively. Constructing the approximate solution is synonymous with ``learning'' a solution and is equally called the ``training'' phase. Similarly, the subsequent use and evaluation of the approximate solution is called the ``inference'' phase. We will follow the common practice of using the training and inference terminology.

It is common for the training and inference stages to be considered as two sequential and independent workloads, whereby a portion of computational time is first spent performing the training only once, and after this is completed the result is subsequently used for inference tasks. When training and inference are overlapped (possibly with feedback between the two), and training can be subsequently revisited, this falls under the categories of reinforcement learning and continual learning (sometimes also called lifelong, incremental, or sequential learning). Situations where the training is done only once and not revisited we will call ``isolated'', following the convention of \citet[\S\,1]{chen2018lifelong}. We will primarily focus on isolated training and not on continual training, but for more background on continual training see \citet{wang2024comprehensive}, \citet{buddhi2023continual}, \citet{hadsell2020embracing}, and \citet{chen2018lifelong}.  

In the context of isolated training, there is a clear divide between the training and inference stages. It is frequently the case that the training can be precomputed ahead of time, either out of convenience or necessity, and such settings are called ``offline'' training. Conversely, if the training computations are performed in the same computational workflow alongside their subsequent use for inference,  then this setting is known as ``online'' training. Offline and online training are both equally commonplace, and it is possible in certain applications to freely use either. If the training phase is computationally expensive compared to the inference, and the subsequent result is intended to be used multiple times in differing settings, then offline training is likely preferable. Conversely, if the training is inexpensive compared to the inference, and the subsequent result is used very few times (possibly only once), then online training is likely preferable. 

In our setting of learning an approximation to the solution of a partial differential equation for subsequent sample path generation, it is easy to envisage situations appropriate for both offline or online training. The isolated offline training setting is the easiest situation to consider, model, and analyse, as the training and inference are two decoupled tasks, which can be handled independently in isolation. Conversely, isolated online training is a more difficult situation to tackle. 

\subsection{Combining multilevel Monte Carlo and neural networks}

For systems of coupled forward-backward stochastic differential equations, it is our ambition to estimate properties of the backward stochastic process, and to be able to sample paths for this process. Using \cref{thm:feynman_kac_semi_linear_multidimensional} and the framework proposed by \citet{raissi2018forwardbackward}, we can use a combination of neural networks and the Euler-Maruyama scheme to sample from the backward process. However, the training and evaluation of neural networks is very expensive, and so care must be taken to avoid any unnecessary calculations using neural networks. The deeper and wider the neural network, the greater the computational cost, and the shallower and narrower the neural network, the cheaper the cost. 

This naturally presents an entry point for applying the multilevel Monte Carlo framework. Small neural networks briefly trained on small amounts of data form crude estimators. Identically, large neural networks thoroughly trained on large amounts of data form the fine estimators. 

A related combination of neural networks and multilevel Monte Carlo methods has been employed by \citet{ko2023multilevel} for gradient estimation of stochastic systems. This though utilised a differing setup to ours, where for the stochastic processes the drift and diffusion functions are themselves neural networks, and hence differs from our setting. This setup is sometimes called a neural stochastic differential equation. Additionally \citet{ko2023multilevel} consider a regular forward-only stochastic differential equation, with no backward stochastic process, and hence a considerably simpler and different setting than our own. Similar work by 
\citet{gierjatowicz2023robust} further investigates such neural stochastic processes, producing efficient Monte Carlo methods for model calibration and training, applying their results in the context of volatility modelling.

As an example of combining differing neural network sophistications into a multilevel Monte Carlo framework, we consider the nested multilevel Monte Carlo \cref{eqt:nested_mlmc}. For a given level of discretisation, let us suppose for simplicity we have only two neural networks available to us: a low fidelity neural network parametrised by some \( \theta \), and a higher fidelity neural network similarly parametrised by some \( \theta' \). What differentiates the low and high fidelities could be e.g.\ the neural networks' sizes, the extent they have been trained, both, etc. What is crucially important is that they are both trained using the same Brownian motion paths, to ensure the two networks are closely coupled, and thus care is needed to properly ensure this (going beyond just seeding random number generators).\footnote{This is indeed one reason to prefer keeping the Brownian motion paths fixed across all training iterations when evaluating the loss function.} Similarly, let us suppose we only entertain using two differing discretisation levels: a fine discretisation, and a coarse one. We denote \( \widehat{P} \) estimated using the fine and coarse discretisations as \( \widehat{P}^{\mathrm{f}} \) and \( \widehat{P}^{\mathrm{c}} \) respectively, and similarly \( \widehat{P}^{\theta} \) and \( \widehat{P}^{\theta'} \) for the differing neural networks, and combine the two superscripts where appropriate.\footnote{The granularity with which the neural network is trained should not be confused with the granularity that the path is being simulated with, and in general these two granularities should be anticipated to be different.} 

We have two obvious means of combining the different levels of discretisation with the differing neural network sophistications. The first is to nest the neural network multilevel setup within the temporal one (or \textit{vice versa}, as both give the same decomposition). The other is to reduce the neural network's fidelity simultaneously with increasing the temporal granularity. These give rise to the multilevel Monte Carlo decompositions 
{\allowdisplaybreaks
\begin{align}
\mathbb{E}(\widehat{P}^{\mathrm{f}, \theta'}) & = 
\begin{aligned}[t]
& \mathbb{E}(\widehat{P}^{\mathrm{c},\theta})
 \\ & {} + \mathbb{E}(\widehat{P}^{\mathrm{c},\theta'} - \widehat{P}^{\mathrm{c},\theta})
\\& {} + \mathbb{E}(\widehat{P}^{\mathrm{f},\theta} - \widehat{P}^{\mathrm{c},\theta})
 \\ & {} + \mathbb{E}(\widehat{P}^{\mathrm{f},\theta'} - \widehat{P}^{\mathrm{c},\theta'} - \widehat{P}^{\mathrm{f},\theta} + \widehat{P}^{\mathrm{c},\theta})
\end{aligned} \label{eqt:nested_mlmc_nn_and_temporal} \\
\shortintertext{and}
\mathbb{E}(\widehat{P}^{\mathrm{f}, \theta'}) & =  \mathbb{E}(\widehat{P}^{\mathrm{c},\theta}) + \mathbb{E}(\widehat{P}^{\mathrm{f},\theta'} -  \widehat{P}^{\mathrm{c},\theta}) \label{eqt:mlmc_nn_and_temporal}
\end{align}
}%
respectively. In \cref{eqt:nested_mlmc_nn_and_temporal} we have ordered the terms in decreasing variance based on our numerical results (shown later in \cref{fig:multilevel_variances}), and we have (for our deferred specific example at least) substantial variance reductions such that 
\begin{multline}
\mathbb{V}(\widehat{P}^{\mathrm{c},\theta})
 \gg \mathbb{V}(\widehat{P}^{\mathrm{c},\theta'} - \widehat{P}^{\mathrm{c},\theta})
 \gg \mathbb{V}(\widehat{P}^{\mathrm{f},\theta} - \widehat{P}^{\mathrm{c},\theta})
\\ \gg \mathbb{V}(\widehat{P}^{\mathrm{f},\theta'} - \widehat{P}^{\mathrm{c},\theta'} - \widehat{P}^{\mathrm{f},\theta} + \widehat{P}^{\mathrm{c},\theta}).
\end{multline}
The behaviour of the variance of the \( \widehat{P}^{\mathrm{c},\theta'} - \widehat{P}^{\mathrm{c},\theta} \) term is the most important from an analysis viewpoint, as it is the leading order multilevel correction term. Our numerical analysis in \cref{sec:numerical_analysis} focuses on producing bounds that we can use for such terms (or proxies thereof). We also gather numerical results for showcasing the empirical behaviour of the variance exhibited for the other higher order correction terms, although exact numerical bounds remain open for further work. 

Empirically we found (but for visual clarity did not show) that 
\begin{equation}
\mathbb{V}(\widehat{P}^{\mathrm{f},\theta'} - \widehat{P}^{\mathrm{c},\theta}) \approx \mathbb{V}(\widehat{P}^{\mathrm{c},\theta'} - \widehat{P}^{\mathrm{c},\theta}), 
\end{equation}
and so we expect \cref{eqt:mlmc_nn_and_temporal} to have an identical performance to \cref{eqt:nested_mlmc_nn_and_temporal}. Thus while the leading order terms in \cref{eqt:nested_mlmc_nn_and_temporal} may be easier to analyse, \cref{eqt:mlmc_nn_and_temporal} should behave identically but be computationally easier to implement. Consequently \cref{eqt:mlmc_nn_and_temporal} is the multilevel Monte Carlo framework we propose.

\subsubsection{Offline and online training}

Isolated offline training allows us to consider the task of learning an approximation to the solution of the underlying partial differential equation as its own independent and individual challenge. This has been the primary setting considered by e.g.\ \citet{raissi2018forwardbackward}, \citet{guler2019robust}, and \citet{e2017deep_learning,e2018solving,e2021algorithms}. Isolated online training forces us to consider the task of learning an approximation to the solution of the underlying partial differential equation alongside its subsequent use for inference. When limited computational resources are available, determining how to optimally \textit{a priori} divide these resources between training and inference is a substantially more difficult setting than the offline training setting. To the best of our knowledge, we have not found any literature appropriate to this setting. 
Continual online training is likely too difficult and remains a topic for further research. 

\subsection{Previous research}
\label{sec:previous_research}

The most relevant previous works we consider are those by \citet{e2017deep_learning,e2018solving,e2021algorithms}, \citet{raissi2018forwardbackward}, and \citet{guler2019robust}, where the former two are most closely aligned to our exact setup.

In the work by \citet{guler2019robust}, their investigations cover two topics. The first is the use of differing neural network architectures in the training phase, reviewing and comparing the stability and generalisability of fully connected conventional deep neural networks, residual neural networks (ResNet) \citep{he2016deep}, and non-autonomous input-output stable networks (NAIS-Net) \citep{ciccone2018nais}. \citet{guler2019robust} empirically demonstrated that NAIS-Net consistently appeared the most favourable in their setup from the viewpoint of stability and generalisability, smoothing out the variability seen during training by \citet{raissi2018forwardbackward} considerably, albeit typically incurring twice the computational cost measured with respect to training time. Interestingly, \citet{guler2019robust} highlight and emphasise that the problem of training a neural network can be viewed as an optimal control problem, as highlighted by \citet{benning2019deep}, \citet{liu2019deep}, and \citet{li2018optimal}.

The second topic covered by \citet{guler2019robust} is the use of a ``multilevel discretisation method'' which is ``inspired'' by the multilevel Monte Carlo discretisation framework introduced by \citet{giles2008multilevel}. We would like to emphasise that their framework is only inspired by the ideas underpinning multilevel Monte Carlo. Unfortunately, neither the paper by \citet{guler2019robust}, nor their underlying code, offer any details of their ``multilevel discretisation method''. From direct conversations with \mbox{P.~Parpas} (a co-author of \citep{guler2019robust}), our understanding is that \citet{guler2019robust} implemented a routine approximately following \cref{algo:mlmc_inspired_training}, as has been subsequently repeated by \citet{naarayan2024thesis}. Consequently, we contribute \cref{algo:mlmc_inspired_training} as our formal specification of what we understand \citet{guler2019robust} to have intended, with one small but key distinction (discussed shortly). 

\begin{algorithm}[htb]
	\DontPrintSemicolon
	\KwIn{Initial prior values \( \theta_\ast \), maximum discretisation level \( L \), and the maximum number of iterations \( K \).}
	\KwOut{Approximate optimiser of \cref{eqt:loss}.}
	\algrule
	\SetKwProg{initialisation}{Initialisation}{}{}
	\initialisation{}{
		Set \( \theta_{-1} \leftarrow \theta_\ast \).\;
		Generate \( M 2^L \) Gaussian random vectors \( \{Z_1, Z_2, \ldots, Z_{M2^L}\} \) each of dimension \( d \).\;
	}{}
	\SetKwProg{iterations}{Training iterations}{}{}
	\iterations{}{
		\For{\( l \in \{0, 1, \ldots, L\} \)}{
			Set \( N \leftarrow 2^{l} \).\;  
			Set \( \Delta t \leftarrow 2^{-l} \).\; 
			For each of the \( 1 < m \leq M \) batches, set \( \Delta W_n^{(m, l)} \leftarrow  2^{-L/2} \sum_{k=1}^{2^{L-l}} Z_{\kappa(n,m,l) + k}\),
			where \( \kappa(n,m,l) \coloneq  2^L (m-1) + n2^{L-l} \). \label{algo:step:construct_coupled_brownian_paths}\;
			Using the starting point \( \theta_{l-1} \), perform \( \tfrac{K}{L+1} \) iterations of an iterative optimisation algorithm to find the minimiser of the loss defined by \cref{eqt:loss}, producing estimator \( \theta_l \).\;
		}
		\Return \( \theta_L \).\;
	}{}
	\caption{A multilevel Monte Carlo inspired training algorithm.}
	\label{algo:mlmc_inspired_training}
\end{algorithm}  

\phantomsection \label{paragraph:loss_function_resampling}
There are a few characteristics of \cref{algo:mlmc_inspired_training} that are worth discussing. The first is that the loss functions optimised between the various levels use the same underlying Brownian motions at each level of iteration. The hope is to ensure that the starting point obtained from the previous level acts as a good starting point for the next. Conversely though, readers more familiar with applied statistics might equally worry about the implications of using the same Brownian motions paths at each level, and the associated risk that comes with over-fitting. Generating new samples at each level could encourage learning a set of model parameters that generalise well, which is what was done in the implementations of  \citet{e2018solving} and also \citet{raissi2018forwardbackward}. This difference is the key distinction between what we propose in \cref{algo:mlmc_inspired_training} and what we believe other authors have implemented. 

The second characteristic of this though is that the model uses the information gained from the coarser level when optimising at the finer level. This both adds a sequential aspect to the entire algorithm, but additionally adds a dependency which breaks the telescoping summation central to the multilevel Monte Carlo method. To make this point clearer, we add a bit more mathematical notation and formalism to pinpoint the violation. Let \( \vartheta \) be the optimal model parameters in general for the loss, where specifically they may be optimal in the sense e.g.\ \( \vartheta \coloneq \lim_{N,M\to\infty} \argmin \mathscr{L} \), where we assume \( \vartheta \) exists and is unique. As \( \vartheta \) is a random variable, we are interested in knowing \( \mathbb{E}(\vartheta) \). This is unfortunately not possible to do exactly, and hence we wish to approximate this by the original loss from \cref{eqt:loss} without the asymptotic limits on the discretisation and batch size. If we suppose that a discretisation using \( N=2^L \) intervals and \( M \) Brownian motion realisations is the finest quality estimator we can afford, then we define \( \theta^{\mathrm{f}} \coloneq \argmin \mathscr{L} \) and we use the Monte Carlo estimator \( \mathbb{E}(\vartheta) \approx \mathbb{E}(\theta^{\mathrm{f}}) \). Similarly, a minimiser estimated using coarser discretisations we would denote \( \theta^{\mathrm{c}} \). As we have seen in \cref{algo:mlmc_inspired_training}, for an iterative scheme, we require an initial value prior estimator for \( \vartheta \). Denoting \( \mathbb{E}(\theta^{\mathrm{f}} \mid \theta_\ast) \) as the estimator produced by \cref{algo:mlmc_inspired_training} given an initial prior value \( \theta_\ast \), we see that we have the estimator \( \mathbb{E}(\vartheta) \approx \mathbb{E}(\theta^{\mathrm{f}} \mid \theta_\ast)  \).
The multilevel Monte Carlo framework to split the estimator would be
{\allowdisplaybreaks
	\begin{subequations}
		\begin{align}
			\mathbb{E}(\vartheta) 
			& \approx \mathbb{E}(\theta^{\mathrm{f}} \mid \theta_\ast) \label{eqt:mc_estimator}\\
			& = \mathbb{E}(\theta^{\mathrm{c}} \mid \theta_\ast) + \mathbb{E}(\theta^{\mathrm{f}} - \theta^{\mathrm{c}} \mid \theta_\ast) \label{eqt:mlmc_estimator} \\
			& \neq \mathbb{E}(\theta^{\mathrm{c}} \mid \theta_\ast) + \mathbb{E}(\theta^{\mathrm{f}} - \theta^{\mathrm{c}} \mid \mathbb{E}(\theta^{\mathrm{c}} \mid \theta_\ast) ) \label{eqt:mlmc_estimator_attempt_with_correction} \\
			& \neq \mathbb{E}(\theta^{\mathrm{f}} \mid \mathbb{E}(\theta^{\mathrm{c}} \mid \theta_\ast)).
			\label{eqt:mlmc_estimator_attempt_no_correction}
		\end{align}
	\end{subequations}
}%
The first approximation in \cref{eqt:mc_estimator} is the  regular Monte Carlo estimation, and the first equality giving \cref{eqt:mlmc_estimator} is the correct way to perform the multilevel Monte Carlo decomposition, where the equality is exactly preserved because we have taken care to ensure we have a telescoping summation. The two inequalities \cref{eqt:mlmc_estimator_attempt_with_correction,eqt:mlmc_estimator_attempt_no_correction} however are not correct multilevel Monte Carlo decompositions, for exactly the reason that they violate the telescoping summation. Here \cref{eqt:mlmc_estimator_attempt_no_correction} is the basis of the estimator produced using \cref{algo:mlmc_inspired_training}, and \cref{eqt:mlmc_estimator_attempt_with_correction} is the attempted multilevel Monte Carlo framework proposed by \citet{guler2019robust} and \citet{naarayan2024thesis}. 

Our intention here is not to discard nor speak disparagingly about the utility of the estimator \cref{eqt:mlmc_estimator_attempt_no_correction} produced using \cref{algo:mlmc_inspired_training}, but only to criticise calling it a multilevel Monte Carlo estimator. In fact, as highlighted by \citet{guler2019robust} and \citet{naarayan2024thesis}, \cref{eqt:mlmc_estimator_attempt_with_correction} provides a good estimator that appears stable and computationally favourable compared to the regular Monte Carlo estimator from \cref{eqt:mc_estimator}. Thus our primary concern is with the apparent misnomer. The theoretical importance and significance of respecting the telescoping sum cannot be understated, as repeatedly emphasised by \citet{giles2015multilevel}, although in practice methods which violate this can still perform well (albeit without having the same theoretical assurances to support them). Consequently, identifying such a misnomer and clarifying the framework used by \citet{guler2019robust} we feel is an important correction and distinction worth emphasizing. Nonetheless, we do empathise with \citet{guler2019robust} and do not have an obvious proposal for a more appropriate name, recognising the difficulty of naming as [in]famously surmised by \mbox{P.~Karlton:} ``\textit{there are only two hard things in computer science: cache invalidation and naming things}''.

\section{Numerical analysis}
\label{sec:numerical_analysis}

In this \namecref{sec:numerical_analysis} are the main contributions of this work, containing a mixture of heuristic and more rigorous analytic treatments, resulting in new and original contributions. We will put the loss function proposed by \citet{e2018solving} and \citet{raissi2018forwardbackward} under the microscope, giving a heuristic discourse of its bias and variance properties. Thereafter, we produce strong error bounds (appropriate for multilevel Monte Carlo) for the estimators produced using the setup by \citet{raissi2018forwardbackward}. These results, while not unsurprising, are necessary to bridge the gap between the empirical findings of \citet{raissi2018forwardbackward} and the analytic assurances numerical analysts (and computing practitioners) strive for, opening the door to comparisons with the wider body of existing frameworks. The necessity to formalise the framework proposed by \citet{raissi2018forwardbackward}, and give convergence theorems to support the empirical results, is well expressed by \citet{lamport2022interview}: ``\textit{When you write an algorithm, you need to have a proof that it's correct. An algorithm without a proof is a conjecture, it's not a theorem}''.

\subsection{Convenient notation}

To reduce notational clutter in the upcoming numerical analysis, we define some convenient shorthands. We begin with those which facilitate concisely indexing various times. 

\begin{definition} For a continuous time interval \( t \in [0, T] \), discretised into the \( N+1 \) points \( 0 = t_0 < t_1 < \cdots < t_{N-1} < t_N = T\), the index of the nearest discretisation time point which is immediately proceeding or equal to a given time \( t \) is denoted 
	\(n_t \coloneq \max \{n\in\{0,1,\ldots,N\} \colon t_n \leq t\}\).
\end{definition}

\begin{definition}
	\label{def:time_subscript_abbreviation}
	For a continuous process \( A_t \) defined for \( t \in [0, T] \), we make the abbreviation \(A_{n} \equiv A_{t_n}\) for the time points \(t_n \in \{0, t_1, t_2, \ldots, t_{N-1}, T\}\). For a discrete process \( B_n \) defined for \( n \in \{0,1,\ldots,N\} \), we make the abbreviation \(B_{t} \equiv B_{n_t}\).
\end{definition}

There can arise points of ambiguity, such as e.g.\ for \( t \in [0,1] \) whether \(\widehat{A}_1\) refers to \(\widehat{A}_{t_1}\) or \(\widehat{A}_{t_N}\). Such ambiguities are easily reconciled by the context, and if this is not the case then we will clarify this at the point of use. 

When we later produce our strong error bounds it will be convenient to extend our approximations from a discrete set of times to a continuous one. 

\begin{definition}
\label{def:continuous_time_interpolation}
For a continuous time process \(A_t\) in the interval \(t \in [0, T] \), approximated by some \(\widehat{A}_{t_n} \approx A_{t_n}\) at the discrete times \(t_n \in \{0, t_1, t_2, \ldots, t_{N-1}, T\}\), the piecewise constant interpolation \citep[\S\,9.1, (1.14)]{kloeden1999numerical} \(\overline{A}_t \coloneq \widehat{A}_{n_t}\)
gives an approximation \(\overline{A}_t \approx A_t \) for any \(t \in [0, T]\).
Similarly, for any function \(f\) which is not an approximation, but whose value at any time \(t_n \in \{0, t_1, t_2, \ldots, t_{N-1}, T\}\) is known and is \(f_{t_n}\), the piecewise constant interpolation is \(\overline{f}_t \coloneq f_{n_t}\) for any \(t \in [0, T]\).
\end{definition}

The interpolation from \cref{def:continuous_time_interpolation} has the nice property that it is mathematically very simple to use, and that it is correctly adapted. However, there is the piecewise linear interpolation \citep[\S\,9.1, (1.16)]{kloeden1999numerical} \(\overline{A}(t) \coloneq \widehat{A}_{n_t} + \tfrac{t-t_{n_t}}{t_{n_t + 1} - t_{n_t}}(\widehat{A}_{n_t + 1} - \widehat{A}_{n_t})\) for \( t \in [t_{n_t}, t_{n_t + 1}] \). This has the advantage of being continuous, but the drawback of not being correctly adapted. Strong error bounds for the continuous time approximation \citep[\S\,10.6]{higham2021introduction} of regular stochastic differential equations using the Euler-Maruyama scheme are slightly weaker than the discrete time error bounds, and the differing bounds for the piecewise constant and piecewise linear approximations are studied by \citet{muller2002optimal}. 

In our upcoming error bounds, we will be interested in the asymptotic behaviour as models become well trained and discretisations become small. In the same manner that we restate \( f \leq 2 \Delta t \) as \( f = O(\Delta t) \) using ``big \( O \)-notation'', we wish to extend this even further, both by dropping the need to clutter the analysis with \( O(\cdot) \) everywhere, but also by absorbing all the less important terms distracting us from the main result we wish to convey.  

\begin{definition}
	\label{def:relational_o_notation}
	For two quantities \( f \) and \( g \) we introduce the relation ``\( \preceq \)''  and write \( f \preceq g \) to represent \( f = O(g) \) where \( g \) can depend only on: the temporal discretisation \( \Delta t \), the neural network's error parameter \( \epsilon_{\theta} \),  the discretisation parameter \( N \), and the iteration index \( n \), but not on the ratio \( N\Delta t \) (which is equal to \( T \)). All other dependencies are subsumed into the coefficients hidden by the implicit \( O(\cdot) \).
\end{definition}

Examples of quantities which are concealed using this notation include (using \( \Delta t \) in these examples): 
\begin{longdescription}
	\item[Integral constants] Trivial integral constants are subsumed, e.g.\ \( f \leq 2 \Delta t \preceq \Delta t\).
	\item[Unbalanced constants] Constants for differing terms are subsumed,
	e.g.\ \( f \leq 2 \Delta t + 3 \Delta t^2 \preceq \Delta t + \Delta t^2 \).
	\item[Continuity constants] For a function \( f \) with Lipschitz constant \( L \) we write \( \lvert f(t + \Delta t) - f(t) \rvert \leq L \Delta t \preceq \Delta t \), and similarly so for the equivalent constant appearing in expressions of H\"{o}lder continuity. 
	\item[Moment constants] By Jensen's inequality we have 
	\( \lvert \sum_{k=1}^{K} f_k \rvert^p  \leq K^{p-1}  \sum_{k=1}^{K} \lvert f_k \rvert^p  \preceq  \sum_{k=1}^{K} \lvert f_k \rvert^p\). This also subsumes the similar terms arising in Doob's inequality and the Burkholder-Davis-Gundy inequality.
	\item[Parametrised constants] The constant \( T \) is subsumed where
	\( f \leq h(T) g \preceq g \). Note that this does not subsume \( n \) nor \( N \). This is relevant to applications of Gr\"{o}nwall's inequality.
\end{longdescription}

\subsection{The loss function}
\label{sec:loss_function_of_raissi}

When we consider the loss function \cref{eqt:loss} used for training the neural network proposed by \citet{raissi2018forwardbackward}, the motivation for this is that it will force us to learn an approximate solution which is accurate along the entire duration of the whole path, not just at the terminal or initial conditions. Implicit with this setup is the suggestion that as more training iterations are performed, the approximate solution should approach the exact solution. However, is this the case?

The loss function in \cref{eqt:loss} looks to compare the backward process' approximate value obtained from using \(\hat{u}\) and \cref{algo:euler_maruyama_with_nn} with that which would have been obtained by an Euler-Maruyama update using \cref{eqt:fbsde_bsde_em}. Supposing after the \(n\)-th iteration we have values \((\widehat{X}^{\ast}_n, \widehat{Y}^{\ast}_n, \widehat{Z}^{\ast}_n)\) where \( \widehat{X}_n^{\ast} \) can denote either \( \widehat{X}_n \) or \( \widehat{X}_n^{\theta} \) and similarly for \( \widehat{Y}^{\ast}_n \) and  \( \widehat{Z}^{\ast}_n \), and we compare the values obtained for the backward process at iteration \(n+1\) from the Euler-Maruyama scheme using \cref{eqt:fbsde_bsde_em}, denoted as \( \widehat{Y}^{\ast,\mathrm{EM}}_n \), against the neural network approximation \(\hat{u}\). The difference between these two estimates using Taylor's theorem is 
\begin{equation}
\label{eqt:difference_em_nn}
\widehat{Y}_{n+1}^{\ast,\mathrm{EM}} - \hat{u}(t_{n+1}, \widehat{X}_{n+1}^{\ast};\theta) = \sum_{i=1}^{9} R_i,
\end{equation}
where \(f^\ast_n\) denotes \(f\)  being evaluated at \( (t_{n}, \widehat{X}^{\ast}_n, \widehat{Y}^{\ast}_n, \widehat{Z}^{\ast}_n) \). Dropping for notational simplicity the explicit transposes, trace operators, and element-wise summations (which should be clear from context), the remainder terms are
{\allowdisplaybreaks
\begin{subequations}
\begin{align}
R_1 & \coloneq \widehat{Y}_n^{\ast} - \hat{u}^{\ast}_n, \\
R_2 & \coloneq (\widehat{Z}_n^{\ast} - b_n\diffp{\hat{u}_n^{\ast}}{x}) \Delta W_n, \\
R_3 & \coloneq  (\phi_n^{\ast} - \diffp{\hat{u}_n^{\ast}}{t} 
- a^{\ast}\diffp{\hat{u}_n^{\ast}}{x} 
- \frac{1}{2}(b_n^{\ast})^2 \diffp[2]{\hat{u}^{\ast}}{x} )\Delta t, \\
R_4 & \coloneq - \frac{1}{2} (b_n^{\ast})^2 \diffp[2]{\hat{u}_n^{\ast}}{x} (\Delta W_n^2 - \Delta t), \\ 
R_5 & \coloneq - \frac{1}{2} (a_n^{\ast})^2 \diffp[2]{\hat{u}_n^{\ast}}{x} \Delta t^2, \\
R_6 & \coloneq - a_n^{\ast} b_n^{\ast} \diffp[2]{\hat{u}_n^{\ast}}{x} \Delta W_n \Delta t, \\
R_7 & \coloneq - \frac{1}{2} \diffp*[2]{\hat{u}(\tau_n, \xi_n)}{t} \Delta t^2, \\
R_8 & \coloneq - \diffp*[2]{\hat{u}(t_n, \xi_n')}{t,x} \Delta t \Delta \widehat{X}_n^{\ast}, \\
\shortintertext{and}
R_9 & \coloneq - \diffp*[3]{\hat{u}(t_n, \xi_n'')}{x} \Delta (\widehat{X}_n^{\ast})^3,
\end{align}
\end{subequations}
}%
for some \(\tau_n \in (t_n, t_{n+1}) \) and \(\xi_n, \xi'_n, \xi''_n \in (\widehat{X}_n^{\ast}, \widehat{X}_{n+1}^{\ast})\) with \(\Delta \widehat{X}_n^\ast \coloneq a_n^{\ast} \Delta t + b^{\ast}_n \Delta W_n\).

\citet{raissi2018forwardbackward} posits that by minimising differences of the form in \cref{eqt:difference_em_nn}, that this is a means of learning the approximation \(\hat{u} \approx u \). Indeed, we can see that this is partially true, whereby trivially \cref{eqt:difference_em_nn} is minimised in the \(\ell^2\)-norm if we can achieve all the \(R_i = 0\), and noting that \(R_2 = R_3 = 0\) specifically results in exactly learning \cref{eqt:feynman_kac_semi_linear_pde} and thus finding the exact partial differential equation solution in \cref{thm:feynman_kac_semi_linear_multidimensional}.

Consequently, we can observe that in fact we only wish to learn the minimiser of \(R_2 + R_3\), but instead are learning the minimiser of the sum including all the other remainders. With the inclusion of all the other terms, we can see that we are including terms we are not concerned with minimising. Thus \cref{eqt:difference_em_nn} is a biased loss function for our learning objective (namely finding the exact solution to \cref{eqt:loss}). Now that we can appreciate the bias of the loss function, we can quantify it. We can first note that of all the terms in \(\{R_i\}_{i=4}^9\), that \(R_4\) is the smallest in expectation. It is a well known result that \( W_t^2 - t \) is a martingale \citep[theorem 3.7]{klebaner2012introduction}, meaning \(R_4\) has zero expectation, and it is straightforward to show \(\mathbb{V}(\Delta W_n^2) = 2 \Delta t^2\), which implies \(\mathbb{E}(\lvert \Delta W^2_n - \Delta t  \rvert) \leq \sqrt{2} \Delta t\) by an application of the Cauchy-Schwarz inequality. Consequently, we can see that \(\mathbb{E}(R_4) = 0\) and \(\mathbb{E}(\lvert R_4 \rvert) = O(\Delta t)\). Similarly \(\mathbb{E}(R_5) \neq 0 \) with \(\mathbb{E}(\lvert R_5 \rvert) = O(\Delta t^2)\). Thus we have a combination of systematic bias terms and unsystematic bias terms. Heuristically by expanding the last few \( \Delta \widehat{X}_n^{\ast} \) terms in the remainders we expect the leading order systematic term to be \(O(\lvert \Delta t^{3/2} \rvert)\) and the bias to fluctuate with variations of size \(O(\lvert \Delta t^{1/2} \rvert)\). However, we can readily see that at the \(n\)-th iteration \(R_4\), \(R_5\), and \(R_6\) are known quantities. 

Readers familiar with the numerical solution of stochastic differential equations may recognise the \( (\Delta W_n^2 - \Delta t) \) term appearing in \( R_4 \), as a similar term appears in the Milstein scheme. The Milstein scheme \citep{milstein1974approximation,milstein1995numerical} is the next higher order numerical scheme beyond the Euler-Maruyama scheme. A derivation can be found in \citet[\S\,10.3]{kloeden1999numerical}, and in 1-dimension when the forward process is decoupled from the other processes produces the well known approximation 
\begin{equation}
\label{eqt:milstein_scheme_1d}
\widehat{X}_{n+1} = \begin{multlined}[t] 
\widehat{X}_n + a(t_n,\widehat{X}_n) \Delta t_n + b(t_n, \widehat{X}_n) \Delta W_n \\
+ \frac{1}{2} b(t_n, \widehat{X}_n) \nabla b(t_n, \widehat{X}_n) (\Delta W_n^2 - \Delta t).
\end{multlined}
\end{equation}
Comparing this to the Euler-Maruyama scheme from \cref{eqt:fbsde_fsde_em} we can see the new final \( (\Delta W_n^2 - \Delta t) \) term. We can compute what might be the equivalent coefficient for the \( (\Delta W_n^2 - \Delta t) \) term for a Milstein scheme approximation to \( Y \), and it is straightforward to show that this is not equal to the same coefficient appearing in \( R_4 \). The implication of this inequality is that the next leading order corrections that would be learnt by minimising the loss are not the corrections introduced by the Milstein scheme. (If \( \widehat{X} \) were produced using the Milstein scheme rather than the Euler-Maruyama scheme then---ignoring the complications of the Milstein Scheme in high dimensions---one might have reason to hope for such an equality).

The consequence of this is that if we needed to learn the solution \( u \) on a finer grid and reduce the discretisation error, then rather than training on a finer grid, one could instead minimise a different loss function that looked to address the remainder \( R_4 \) (and possibly also \( R_5 \), and \( R_6 \)). The gradient \( \nabla u \) is computed using automatic differentiation,\footnote{For an introduction and survey of automatic differentiation, we recommend \citet{baydin2018automatic}.} and so if one is willing to expend the extra computational effort to similarly compute the Hessian \( \nabla^2 u \) (which for modest dimensionalities might not be prohibitive),\footnote{In \cref{eqt:loss_higher_order} we see that we are eventually computing Hessian-vector products, which are computational feasible, as highlighted by \citet[p.\,17]{baydin2018automatic} and \citet{dixon1991use}.} then if we define some other hidden process
{\allowdisplaybreaks
\begin{subequations}
\begin{align}
H_t & \coloneq \nabla^2 u(t, X_t) \\
\shortintertext{and similarly}
\widehat{H}^{\theta}_n & \coloneq \nabla^2 \hat{u}(t_n, \widehat{X}_n^{\theta};\theta)
\end{align}
\end{subequations}
}%
then we could construct the loss function (written in 1-dimensional form for convenience) 
\begin{align*}
\mathscr{L} & \coloneq  
\begin{multlined}[t]
\sum_{m=1}^{M} \sum_{n=0}^{N-1} \lvert \widehat{Y}^{(m)}_{n+1} - \widehat{Y}^{(m)}_{n} - \Phi^{(m)}_{n}\Delta t \\
\begin{aligned}[t]
& {}- \widehat{Z}^{(m)}_{n} \Delta W^{(m)}_{n} \\
& {}- \frac{1}{2}(b^{\ast}_n)^2 \widehat{H}^{(m)}_{n} ((\Delta W^{(m)}_{n})^2 - \Delta t) \rvert^2 
\end{aligned}
\end{multlined}
 \\
& \phantom{\coloneq} {} + \sum_{m=1}^{M} \lvert \widehat{Y}_N^{(m)} - g(\widehat{X}_N^{(m)}) \rvert^2.  \tageq \label{eqt:loss_higher_order}
\end{align*}
 
So what is the advantage of the loss from \cref{eqt:loss_higher_order} compared to \cref{eqt:loss}? A comparison of the two different loss functions is given in \cref{tab:comparison_of_loss_functions}, whereby we can see that while both have the same bias, the size of \cref{eqt:loss_higher_order} is lower than \cref{eqt:loss}, so thus presents a lower variance objective function to minimise. The results in \cref{tab:comparison_of_loss_functions} would stay the same if the \( R_5 \) and \( R_6 \) terms were included in \cref{eqt:loss_higher_order}, and thus they are omitted from \cref{eqt:loss_higher_order}. By considering the \( R_5 \) term, the bias could reasonably be anticipated to persist even across different batches of Brownian path samples at each training iteration. 

\begin{table}[htb]
\centering
\begin{tabular}{ccc}
Loss & \( \mathbb{E}(\cdot)  \) & \( \mathbb{E}(\lvert \cdot \rvert) \) \\
\hline
\firstrowspacing
\cref{eqt:loss} & \( O(\Delta t^{3/2}) \) & \( O(\Delta t) \) \\
\cref{eqt:loss_higher_order} & \( O(\Delta t^{3/2}) \) & \( O(\Delta t^{3/2}) \) 
\end{tabular}
\caption{Comparisons of loss functions.}
\label{tab:comparison_of_loss_functions}
\end{table}

A last remark we can make here is that \cref{tab:comparison_of_loss_functions} suggests that there is a convergence \( \hat{u} \to u \) with \cref{eqt:loss} at a rate \( O(\Delta t) \). However, as mentioned, for the Euler-Maruyama scheme one would expect a convergence rate of \( O(\Delta t^{1/2}) \). To further assess any rates of convergence between \( \hat{u} \) and \( u \) would require considerably less variable neural network architectures, such as those explored by \citet{guler2019robust} (e.g.\ NAIS-Net).

Just as the previous analysis inspected (using Taylor expansions) how minimising the loss function should cause the solution to be learnt, we can approach this same query from the other direction (using Ito-Taylor expansions), and assess the loss which would manifest even from a perfectly learnt solution. 

Suppose we had access to the exact solution of the forward process \( X_t \) and also the hidden process \( Z_t \), not forgetting from \cref{thm:feynman_kac_semi_linear_multidimensional} that \(Z_t\) depends on \( u \) as specified by \cref{eqt:data_dependence_on_pde_solution}. We can define the local error counterpart of \cref{eqt:loss} arising from using the Euler-Maruyama scheme at the \(n\)-th time point as \(E_n\), where dropping the path realisation superscript for brevity, gives
\begin{equation}
E_{n+1} 
\coloneq  
\int_{t_n}^{t_{n+1}} \dl{{u(s, X_s)}} 
 - \phi_{n}\Delta t 
 - Z_{n} \Delta W_{n}.
\end{equation}
Following the notation of \citet[\S\,5]{kloeden1999numerical}, we introduce the operators 
{\allowdisplaybreaks
\begin{subequations}
\begin{align}
L_0 &\coloneq  \diffp{}{t} + a \diffp{}{x} + \frac{1}{2}b^2\diffp[2]{}{x}, \\
\shortintertext{and}
L_1 & \coloneq  b\diffp{}{x},
\end{align}
\end{subequations}
}%
and so performing an Ito-Taylor expansion \citep[\S\,5]{kloeden1999numerical} of \(u\) gives
\begin{equation}
E_{n+1} = 
\begin{aligned}[t]
& (L_0 u(t_n, {X}_{t_n}) - \phi_{n})\Delta t \\
& {} + (L_1 u(t_n, {X}_{t_n}) - Z_{n}) \Delta W_{n} \\
& {} +   \int_{t_n}^{t_{n+1}} \int_{t_n}^{s} L_0 L_0 u(z, {X}_z)\dl{z} \dl{s} \\
& {} +  \int_{t_n}^{t_{n+1}} \int_{t_n}^{s} L_1 L_0 u(z, {X}_z)\dl{W_z} \dl{s} \\
& {} +   \int_{t_n}^{t_{n+1}} \int_{t_n}^{s} L_0 L_1 u(z, {X}_z)\dl{z} \dl{W_s} \\
& {} +   \int_{t_n}^{t_{n+1}} \int_{t_n}^{s} L_1 L_1 u(z, {X}_z)\dl{W_z} \dl{W_s}.
\end{aligned}
\end{equation}
If the solution to the underlying partial differential equation is perfectly learnt, then by definition, from \cref{thm:feynman_kac_semi_linear_multidimensional} the leading order \(\Delta t\) and \(\Delta W_n\) coefficient terms will be exactly zero. This means that the integral terms will be the only non-zero terms. The integral terms are themselves still difficult to tackle, but are typically considered to be lower order than the leading order \(\Delta t\) and \(\Delta W_n\) terms, with the exception that the final integral might be of the same order as the \(\Delta t\) term (this observation is the basis of producing the Milstein scheme). Consequently, if the solution to the partial differential equation is not perfectly learnt, then we can precisely and separately see the learning errors and the  discretisation errors.

\subsection{Strong error bounds}

When we consider the loss defined by \cref{eqt:loss}, we can see that this corresponds to the \(\ell^2\)-norm of the error, defined as the difference between the estimates of the backward stochastic process produced by the learnt approximate solution of the partial differential equation and that produced by the Euler-Maruyama scheme. For readers more familiar with numerical solutions of stochastic differential equations, there is considerable literature surrounding bounding the expected values of such \(L^2\) errors. Specifically, such errors are measured using the strong error, which is central to multilevel Monte Carlo analyses. This crucial nature of the strong convergence in multilevel Monte Carlo is repeatedly emphasised by \citet{kloeden2013convergence}. The strong error can typically be measured using either the \(L^1\), \(L^2\), or \(L^p\)-norms, and is either evaluated at the terminal value, or via the supremum over intermediate values. Thus for a stochastic process \(A_t\) for times \(t \in [0,T]\) with some approximation \(\widehat{A}_{n}\) where \(\widehat{A}_n \approx A_{t_n} \) for the intermediate times \( 0=t_0 < t_1 < \cdots < t_N=T\), the strong error is typically measured using either 
{\allowdisplaybreaks
\begin{subequations}\label{eqt:strong_errors}
\begin{gather}
\mathbb{E}(\lvert A_T - \widehat{A}_N \rvert ),	 \label{eqt:strong_error_weakest} \\
\mathbb{E}(\lvert A_T - \widehat{A}_N \rvert^2 ), \\
\sup_{n\leq N}\mathbb{E}(\lvert A_{n} - \widehat{A}_n \rvert^2 ), \\
\mathbb{E}(\sup_{n\leq N}\lvert A_{n} - \widehat{A}_n \rvert^2 ),  \label{eqt:strong_error}\\
\shortintertext{or}
\mathbb{E}(\sup_{n\leq N}\lvert A_{n} - \widehat{A}_n \rvert^p ) \label{eqt:strong_error_lp}
\end{gather}
\end{subequations}
}%
for any \(p \geq 1\).\footnote{For multi-dimensional processes swap \(\lvert \cdot \rvert^p \to \lVert \cdot \rVert^p_p\).} If the strong error is bounded from above by a term proportional to \(\Delta t^\gamma\) for some \(\gamma > 0\) in the limit \(\Delta t \to 0\), then we say the approximation converges strongly with order \(\gamma\). 

The strong errors presented in \cref{eqt:strong_errors} are in the order of the weakest to the strongest forms.\footnote{We can bound an \(L^q\)-norm by an \(L^p\)-norm for \(1\leq q \leq p\) by use of H\"{o}lder's inequality. Similarly, strong errors using the supremum inside the expectation can be bound by those without it by using Doob's inequality.} For our purposes, we will regard them as interchangeable for most practical purposes, and speculatively [and optimistically] assume that bounds on one translate to equivalent bounds on another with no (or minimal) substantive changes. Such an approach and alternating use of differing strong error definitions is not unusual (cf.\ \citep[p.\,85]{higham2021introduction}), and often the results do carry over between differing definitions, (e.g.\ \citep[theorem 10.2.2]{kloeden1999numerical}). Similarly \(\ell^p\) bounds regularly match closely \(L^p\) bounds (cf.\ \citep[\S\,10.6]{higham2021introduction}) when extended through either piecewise constant or piecewise linear interpolations \citep{muller2002optimal}. 

The preeminent analysis bounding the strong error for the Euler-Maruyama scheme is by \citet[\S\,10]{kloeden1999numerical}, although a more accessible introduction is given by \citet[\S\,10]{higham2021introduction}. A generalised extension to this framework, targeting altered Euler-Maruyama schemes is also presented by \citet[lemma~4.3]{giles2022approximate}.\footnote{ \citeauthor{giles2022approximate} use this extension to provide strong error bounds for altered Euler-Maruyama schemes using approximations for random variable distributions \citep{giles2022approximate,giles2023approximating} or those incorporating numeric rounding error \citep{giles2024rounding}.} The proofs of these frameworks proceed by the same steps: performing Ito-Taylor expansions, utilising Lipschitz continuity of the drift and diffusion functions, and lastly using a combination of the Burkholder-Davis-Gundy inequality \citep{burkholder1972integral}, Doob's inequality, and Gr\"{o}nwall's inequality \citep{gronwall1919derivatives}.  

\subsubsection{Strong error bounds for decoupled processes}

Ideally we would like to understand the strong convergence of our \( \widehat{Y} \) process when the various processes are coupled together (hoping that the approximation scheme does indeed converge at all). However, many systems which are either decoupled or forward-only naturally arise in applications, and consequently there is an appreciable collection of associated analytic results. 

For convergence rates of multiple coupled forward-only stochastic processes, there are the works by \citet{cozma2018strong} and \citet{cozma2018convergence}. However, these and related discussions \citep{kloeden2013convergence,cozma2018strong,cozma2018convergence,higham2002strong,higham2005convergence,giles2009analysing} typically focus on non-linear drift and diffusion processes, and where the emphasis is primarily concerning violations of Lipschitz continuity, rather than on the difficulties of coupling forward and backward processes. 

A recent analysis concerning the Heston model which utilises neural networks is provided by \citet{perotti2024pricing}, which itself builds on the related work of \citet{yarotsky2017error}. While the application does not utilise neural networks for producing approximations to any backward processes, it does give a flavour of the types of interesting results which arise from incorporating neural networks. A fascinating result of this work are bounds on the neural network's size and depth for convergence \citep[theorem 4.6]{perotti2024pricing}. We restate their result, and wish to emphasize the utility for practitioners of having expressions for the neural network's size and depth.

\begin{theorem}
\label{thm:perotti2024pricing}
(\citeauthor{perotti2024pricing} (\citeyear{perotti2024pricing}) \citep{perotti2024pricing})
Let \( \mathscr{F} \) be the space of functions which are \( C^{a-1}([0,1]^b) \) for any choice of \( a,b \in \mathbb{N}\backslash\{0\} \) and whose derivatives up to the \( (a-1) \)-th order are Lipschitz continuous and equipped with the norm defined in \citep[(29)]{perotti2024pricing}. Then for any \( \epsilon \in (0,1) \) there exists a neural network architecture \( H \) using ReLU activation functions such that \( H \) can approximate any \( f \in \mathscr{F} \) with an error less than \( \epsilon \), and \( H \) has at most \( c(1-\log(\epsilon)) \) layers and at most \( c \epsilon^{-b/a}(1-\log(\epsilon)) \) neurons for an appropriate constant \( c \) which depends on only \( a \) and \( b \). 
\end{theorem}

Earlier work also by \citet{perotti2022fast_sampling} also similarly focuses on applications utilising deep learning, investigating sampling techniques for time-integrated stochastic bridges, where the start and end conditions of the underlying process are known. 

It is well known \citep[\S\,10]{kloeden1999numerical} that for regular stochastic processes, under appropriate assumptions, the Euler-Maruyama scheme converges strongly with order \(\tfrac{1}{2}\), and similarly the Milstein scheme with order 1. Related error bounds for backward stochastic differential equations are provided by \citet{bouchard2009strong}. It is natural to wonder then if our approximation for the backward stochastic process similarly demonstrates any strong convergence, or if such a result can be readily obtained analytically? While \citet{raissi2018forwardbackward} explored the setup empirically, no supporting analysis was offered, and thus our results fill this important gap. For the simple case where the forward process is decoupled from the rest, we present \cref{thm:strong_error_decoupled_forward} to showcase that the usual Euler-Maruyama strong convergence order carries over to the backward process.

\begin{theorem} 
\label{thm:strong_error_decoupled_forward}
If the forward process \(X_t\) is decoupled from both the backward process \(Y_t\) and hidden process \(Z_t\) such that \(a\) and \(b\) have no dependence on neither \(Y_t\) nor \(Z_t\), then for \(\widehat{Y}_n\) defined by \cref{algo:euler_maruyama_with_nn} we have the bound
\begin{equation}
\mathbb{E}( \sup_{n\leq N} \lvert Y_{n} - \widehat{Y}_n \rvert^p) \preceq \Delta t^{p/2}.
\end{equation}
\end{theorem}

\begin{proof} 
We can directly use \cref{eqt:backward_process_forward_value_definition} to bound the quantity
\begin{align}
\mathrlap{\mathbb{E}(\sup_{n\leq N} \lvert Y_{n} - \widehat{Y}_n \rvert^p)} 
\hspace{2em} & \notag \\
& = \mathbb{E}(\sup_{n\leq N} \lvert u(t_n, X_n) - u(t_n, \widehat{X}_n) \rvert^p) \\
\shortintertext{using the mean value theorem}
& = \mathbb{E}(\sup_{n\leq N} \lvert (\nabla u(t_n, \xi_n))^{\top} (X_n - \widehat{X}_n) \rvert^p) \\
& \preceq \begin{multlined}[t] \mathbb{E}(\sup_{n\leq N} \lVert \nabla u(t_n, \xi_n)\rVert^p_p \\ \times \sup_{n\leq N} \lVert X_n - \widehat{X}_n \rVert^p_p)
\end{multlined}
 \\
\shortintertext{using the Cauchy-Schwarz inequality}
& \preceq \begin{multlined}[t]
(\mathbb{E}(\lvert\sup_{n\leq N} \lVert \nabla u(t_n, \xi_n)\rVert^p_p\rvert^2))^{1/2} \\ \times (\mathbb{E}(\lvert\sup_{n\leq N} \lVert X_n - \widehat{X}_n \rVert^p_p \rvert^2))^{1/2}
\end{multlined} \\
& \preceq \begin{multlined}[t]
(\mathbb{E}(\sup_{n\leq N} \lVert \nabla u(t_n, \xi_n)\rVert^{2p}_p))^{1/2} \\ \times (\mathbb{E}(\sup_{n\leq N} \lVert X_n - \widehat{X}_n \rVert^{2p}_p))^{1/2}
\end{multlined} \\
\shortintertext{using \cref{eqt:pde_solution_lipschitz_bounds}}
& \preceq \begin{multlined}[t]
(\mathbb{E}(\sup_{n\leq N} (1 + \lVert \xi_n\rVert^{2p}_p)))^{1/2} \\ \times (\mathbb{E}(\sup_{n\leq N} \lVert X_n - \widehat{X}_n \rVert^{2p}_p))^{1/2}
\end{multlined} \\
\shortintertext{as \( X_n \) and \( \widehat{X}_n \) have bounded moments \citep[theorem 4.5.4]{kloeden1999numerical}}
& \preceq (\mathbb{E}(\sup_{n\leq N} \lVert X_n - \widehat{X}_n \rVert^{2p}_p))^{1/2} \\
\shortintertext{using the standard strong convergence order of \( \tfrac{1}{2} \) \citep[theorem 10.2.2]{kloeden1999numerical}}
& \preceq \Delta t^{p/2}, 
\end{align}
which completes the proof. \qedhere
\end{proof}

The result from \cref{thm:strong_error_decoupled_forward} is reassuring, but it makes the very strong assumption of the forward processing being decoupled from the backward process, which is not ideal, and something we would like to relax. However, why is \cref{thm:strong_error_decoupled_forward} worth mentioning at all, and is it really a new result? It is a new result because the backward process approximation is defined using \cref{algo:euler_maruyama_with_nn}, and worth mentioning because it is a crucial (albeit unsurprising) result required to bridge the gap between analytic results and the empirical findings of \citet{raissi2018forwardbackward}. The structural form of the result though is unsurprising, (which is in itself reassuring), as it resembles well known similar results \citep[theorem 10.2.2]{kloeden1999numerical}.

Considering the Black-Scholes-Barenblatt equation \cref{eqt:black_scholes_barenblatt_pde}, because it has a closed form solution \cref{eqt:black_scholes_barenblatt_pde_exact_solution} to compare against, we can empirically assess this using \cref{fig:multilevel_variances} (\(\Circle\)~markers), although we defer the discussion until \cref{sec:numerical_results}.

\subsubsection{Strong error bounds for coupled processes}

While coupled processes are intrinsically more difficult then uncoupled processes, their is no shortage of research activity aimed at producing strong error bounds, including when neural networks are utilised. 
\citet{han2020convergence} consider the loss at the terminal time point, and try to learn drift and diffusion proxies for the backward process, (and thus their method of adopting neural networks differs from our setup). Similar setups are also used by \citet{andersson2023convergence}, and most recently by \citet{negyesi2024generalizedconvergencedeepbsde}, where \citet[theorem 2]{negyesi2024generalizedconvergencedeepbsde} produce the usual strong convergence order \( \tfrac{1}{2} \) result for each of the \( \widehat{X} \), \( \widehat{Y} \), and \( \widehat{Z} \) processes.  We restate their result, as it will help to put our results into a broader context for comparison. (\citet{reisinger2023posteriori} have also recently produced $L^2$-error estimators for fully coupled McKean-Vlasov forward-backward stochastic differential equations, whose results predate \citet{negyesi2024generalizedconvergencedeepbsde}, and are in most parts more general).

\begin{theorem}
\label{thm:negyesi2024generalizedconvergencedeepbsde}
(\citeauthor{negyesi2024generalizedconvergencedeepbsde} (\citeyear{negyesi2024generalizedconvergencedeepbsde}) \citep{negyesi2024generalizedconvergencedeepbsde})
For a sufficiently small \( \Delta t \), the neural network discretisations \( (\widehat{X}^{\theta}, \widehat{Y}^{\theta}, \widehat{Z}^{\theta}) \) defined in \citep[(4)]{negyesi2024generalizedconvergencedeepbsde} are \( L^2 \) and satisfy
\begin{multline} 
\sup_{t \in [0, T]} (\mathbb{E}(\lVert X_t - \overline{X}^{\theta}_t \rVert^2) + \mathbb{E}(\lVert Y_t - \overline{Y}^{\theta}_t \rVert^2))  
 \\ + \int_0^T \mathbb{E}(\lVert Z_t - \overline{Z}^{\theta}_t \rVert^2)  
\preceq \Delta t.
\end{multline}
\end{theorem}

With the results of \citet{negyesi2024generalizedconvergencedeepbsde} restated, we present our own analogous strong error bound, appropriate for our setup. For notational convenience, we restrict our attention to the one dimensional case, (removing the clutter of norms and transposes). 

\begin{definition}
The \textit{standard assumptions} are that \( a \), \( b \), and also \( \phi \) each satisfy  \crefrange{asmp:measurability}{asmp:convergent_initial_value}, that \cref{asmp:invertible_and_bounded_diffusion,asmp:bsde_standard_lipschitz_data} are satisfied, that \( u \) satisfies \cref{eqt:pde_solution_lipschitz_bounds}, and furthermore that \( u \) satisfies the bound\footnote{The bound \cref{eqt:pde_solution_lipschitz_bounds_no_bounded_growth} is stricter than \cref{eqt:pde_solution_lipschitz_bounds}.}
\begin{equation}
\label{eqt:pde_solution_lipschitz_bounds_no_bounded_growth}
\sup_{(t,x) \in [0, T]\times \mathbb{R}^d} (\lvert u \rvert^2 + \lVert \nabla u \rVert^2_2)  \leq K
\end{equation} 
for some positive constant \( K \).
\end{definition}

\begin{theorem}
\label{thm:strong_error_exact_em}
For approximations produced using \cref{algo:euler_maruyama_with_nn}, 
under the standard assumptions and the assumptions that \( d = 1 \) and 
\( \mathbb{E}
(\lvert \widehat{X}_{0} - X_{0} \rvert^2 + \lvert \widehat{Y}_{0} - Y_{0} \rvert^2)  = 0\), then
\begin{equation}
\mathbb{E}(\lvert X_n - \widehat{X}_n \rvert^2 + \lvert Y_n - \widehat{Y}_n \rvert^2) \preceq \Delta t .
\end{equation}
\end{theorem}

\begin{proof} 
Following the steps similar to the usual procedure for bounding the strong error of the Euler-Maruyama scheme \citep[\S\,10.3]{higham2021introduction}, for the continuous approximations from \cref{def:continuous_time_interpolation},
it is straightforward to show for a time \(s \in [0, T]\) that
\begin{subequations}
\begin{multline}
\mathrlap{\overline{X}_{s} - X_{s}} 
\\ = 
\begin{aligned}[t]
& \overline{X}_0 - X_0 \\
& +  \int_{0}^{t_{n_s}} (a(\overline{r}, \overline{X}_r, \overline{Y}_r) - a(r, X_r, Y_r)) \dl{r}\\
& +   \int_{0}^{t_{n_s}}( b(\overline{r}, \overline{X}_r, \overline{Y}_r) - b(r, X_r, Y_r)) \dl{W_r} \\
& -   \int_{t_{n_s}}^{s} a(r, X_r, Y_r)  \dl{r} \\
& -   \int_{t_{n_s}}^{s} b(r, X_r, Y_r) \dl{W_r},  
\end{aligned}
\end{multline}
and similarly
\begin{multline}
\mathrlap{\overline{Y}_{s} - Y_{s}} \\ = 
\begin{aligned}[t]
& \overline{Y}_0 - Y_0 \\
& +   \int_{0}^{t_{n_s}} (\phi(\overline{r}, \overline{X}_r, \overline{Y}_r) - \phi(r, X_r, Y_r)) \dl{r} \\
& +   \int_{0}^{t_{n_s}} (\overline{Z}_r - Z_r) \dl{W_r} \\
& -   \int_{t_{n_s}}^{s} \phi(r, X_r, Y_r) \dl{r} \\
& -   \int_{t_{n_s}}^{s} Z_r \dl{W_r},
\end{aligned}
\end{multline}
\end{subequations}
where \( \overline{Z}_r \coloneq \nabla u(\overline{r}, \overline{X}_r) \).
The proof will follow by bounding each of these integrals, and then by defining the related process 
\begin{equation}
D_s \coloneq \mathbb{E}
(\lvert \overline{X}_{s} - X_{s} \rvert^2 + \lvert \overline{Y}_{s} - Y_{s} \rvert^2)
\end{equation}
we will later apply Gr\"{o}nwall's inequality to \(D_s\) to obtain the desired result. 

Bounding the necessary integrals in the \(L^2\)-norm, using the Cauchy-Schwarz inequality we first have
\begin{align}
& \mathbb{E} ( \lvert \int_{0}^{t_{n_s}} (a(\overline{r}, \overline{X}_r, \overline{Y}_r) - a(r, X_r, Y_r)) \dl{r} \rvert^2 ) \notag \\
& \leq t_{n_s} \mathbb{E} (  \int_{0}^{t_{n_s}} \lvert a(\overline{r}, \overline{X}_r, \overline{Y}_r) - a(r, X_r, Y_r) \rvert^2 \dl{r}  ) \\
\shortintertext{using Fubini's theorem}
& \leq t_{n_s}   \int_{0}^{t_{n_s}} \mathbb{E} (\lvert a(\overline{r}, \overline{X}_r, \overline{Y}_r) - a(r, X_r, Y_r) \rvert^2  ) \dl{r} \\
\shortintertext{using \(t_{n_s} < T\) and \cref{asmp:lipschitz_continuity,asmp:temporal_holder_continuity}}
& \preceq \int_{0}^{t_{n_s}} D_r \dl{r}  + \int_{0}^{t_{n_s}} \lvert \overline{r} - r \rvert \dl{r} \\
\shortintertext{using \(t_{n_s} \leq s \leq T \) and \(\lvert \overline{r} - r \rvert \leq \Delta t\)}
& \preceq \int_{0}^{s} D_r \dl{r}  + \Delta t.
\end{align}
An identical bound follows for the integral containing the \(\phi\) difference in its integrand. 

For the Ito integral we can use Ito isometry to produce the bound
\begin{align}
& \mathbb{E} ( \lvert \int_{0}^{t_{n_s}} (b(\overline{r}, \overline{X}_r, \overline{Y}_r) - b(r, X_r, Y_r)) \dl{W_r} \rvert^2 ) \notag \\
& =  \int_{0}^{t_{n_s}} \mathbb{E} (\lvert b(\overline{r}, \overline{X}_r, \overline{Y}_r) - b(r, X_r, Y_r) \rvert^2) \dl{r} \\
\shortintertext{as before}
& \preceq \int_{0}^{s} D_r \dl{r}  + \Delta t.
\end{align}

For the integral of the difference of the hidden process, using Ito isometry we have
\begin{align}
& \mathbb{E} ( \lvert \int_{0}^{t_{n_s}} (\overline{Z}_r - Z_r) \dl{W_r} \rvert^2 ) \notag \\
& = \int_{0}^{t_{n_s}} \mathbb{E}(\lvert \overline{Z}_r - Z_r \rvert^2) \dl{r} \\
\shortintertext{using \cref{eqt:data_dependence_on_pde_solution}}
& = \int_{0}^{t_{n_s}} \mathbb{E}(\lvert \nabla u(\overline{r}, \overline{X}_r) - \nabla u(r, X_r) \rvert^2) \dl{r} \\
\shortintertext{as \(\overline{X}_{r'}\) is a constant for \(r' \in [\overline{r}, r]\), then using the mean value theorem and Jensen's inequality}
& \preceq \begin{multlined}[t]
\int_{0}^{t_{n_s}} \mathbb{E}(\lvert \nabla u(r, \overline{X}_r) - \nabla u(r, X_r) \rvert^2) \dl{r} \\ + \int_{0}^{t_{n_s}}  \int_{\overline{r}}^{r} \mathbb{E}( \lvert \diffp*{\nabla u (r', \overline{X}_{r'})}{t}  \rvert^2) \dl{r'} \dl{r} 
\end{multlined} \\
\shortintertext{under our standard assumptions}
& \preceq \int_{0}^{s} D_r \dl{r} + \Delta t.
\end{align}

By a simple application of the Cauchy-Schwarz inequality we also obtain the bound
\begin{multline}
\mathbb{E} ( \lvert  \int_{t_{n_s}}^{s} a(r, X_r, Y_r) \dl{r} \rvert^2 )
 \\ \leq (s - t_{n_s}) \int_{t_{n_s}}^{s}  \mathbb{E}(\lvert a(r, X_r, Y_r) \rvert^2) \dl{r} \\ \preceq  \Delta t^2 \preceq \Delta t.
\end{multline}
We can obtain an identical bound for all the remaining integrals using a combination of the Cauchy-Schwarz inequality, Ito isometry, and \cref{eqt:pde_solution_lipschitz_bounds}.

Combining all our bounds together, we obtain 
\begin{align}
0 \leq D_s & \preceq  D_0 + \Delta t + \int_{0}^{s} D_r \dl{r}.
\shortintertext{using Gr\"{o}nwall's inequality}
D_s & \preceq D_0 + \Delta t \preceq \Delta t,
\end{align}
completing the proof. \qedhere
\end{proof}

\Cref{thm:strong_error_exact_em} is ideal for showing convergence when we have access to \( u \). However, we should expect \( u \) to be inaccessible in general, and instead for us to only have the approximation \( \hat{u} \). Consequently, we can separately consider the convergence when utilising \( \hat{u} \), and then bootstrap the resultant convergences together (using variants of the triangle inequality). To do so though, we will require certain assumptions about our approximation. 

\begin{assumption}
\label{asmp:lipschitz_approximations}
\( \hat{u} \) is differentiable and satisfies \cref{eqt:pde_solution_lipschitz_bounds_no_bounded_growth}.
\end{assumption}

\begin{assumption}
\label{asmp:uniform_convergence}
For any sufficiently small tolerance \( \varepsilon > 0 \) such that the loss \( \mathscr{L} \) from \cref{eqt:loss} (with the possible inclusion of \cref{eqt:loss_extra_term}) satisfies \( \mathscr{L} \leq \varepsilon \), then \( N \) and \( M \) are sufficiently large so that we have the uniform convergence bound
\begin{equation}
\label{eqt:uniform_convergence_bound}
\sup_{t\in[0, T]} \sup_{x\in\mathbb{R}^d} \lvert u(t, x) - \hat{u}(t, x) \rvert \leq \epsilon_{\theta}
\end{equation}
for some positive constant \( \epsilon_{\theta} \) depending on \( \varepsilon \) and \( N \).
\end{assumption}

\begin{remark}
\Cref{sec:loss_function_of_raissi} discussed at length that the loss function from \cref{eqt:loss} is biased, and thus to assert any notions of uniformity on our approximation we require \cref{asmp:uniform_convergence}. 
\end{remark}

\begin{theorem}
\label{thm:strong_convergence_between_em_and_nn}
Under the same assumptions as \cref{thm:strong_error_exact_em} and \cref{asmp:lipschitz_approximations,asmp:uniform_convergence} 
\begin{equation}
\mathbb{E}(\lvert \widehat{X}_n - \widehat{X}_n^{\theta} \rvert^2 + \lvert \widehat{Y}_n - \widehat{Y}_n^{\theta} \rvert^2) \preceq \epsilon_{\theta}^2.
\end{equation}
\end{theorem}

\begin{proof}
We begin by considering the difference for the backward process approximations
\begin{align}
\widehat{Y}_n - \widehat{Y}_n^{\theta} & \coloneq u(t_n, \widehat{X}_n) - \hat{u}(t_n, \widehat{X}_n^{\theta}; \theta) \\
\shortintertext{introducing a telescoping difference}
& = \begin{aligned}[t]
& u(t_n, \widehat{X}_n) - \hat{u}(t_n, \widehat{X}_n; \theta) \\
& {} + \hat{u}(t_n, \widehat{X}_n;\theta) - \hat{u}(t_n, \widehat{X}_n^{\theta}; \theta)
\end{aligned}
\\
\shortintertext{using Taylor's theorem}
& = 
\begin{aligned}[t]
& u(t_n, \widehat{X}_n) - \hat{u}(t_n, \widehat{X}_n; \theta) \\
& {} - \nabla \hat{u}(t_n, \xi_n;\theta) (\widehat{X}_n - \widehat{X}_n^{\theta})
\end{aligned} 
\end{align}
using \cref{asmp:lipschitz_approximations,asmp:uniform_convergence}
\begin{equation}
\lvert \widehat{Y}_n - \widehat{Y}_n^{\theta} \rvert \preceq \epsilon_{\theta}  + \lvert \widehat{X}_n - \widehat{X}_n^{\theta} \rvert. \label{eqt:backward_learning_diff_bound}
\end{equation}
Similarly, we can obtain the relation
\begin{multline}
\widehat{X}_{n+1} - \widehat{X}_{n+1}^{\theta} 
\\ = \begin{aligned}[t]
& \widehat{X}_{n} - \widehat{X}_{n}^{\theta} \\
& + \diffp*{a(t_n, \xi'_n, \widehat{Y}_n)}{x} (\widehat{X}_{n}^{\theta} - \widehat{X}_{n}) \Delta t \\
& + \diffp*{a(t_n, \widehat{X}_n^{\theta}, \zeta'_n)}{x} (\widehat{Y}_{n}^{\theta} - \widehat{Y}_{n}) \Delta t \\
& + \diffp*{b(t_n, \xi''_n, \widehat{Y}_n)}{x} (\widehat{X}_{n}^{\theta} - \widehat{X}_{n}) \Delta W_n \\
& + \diffp*{b(t_n, \widehat{X}_n^{\theta}, \zeta''_n)}{x} (\widehat{Y}_{n}^{\theta} - \widehat{Y}_{n}) \Delta W_n.
\end{aligned} \label{eqt:forward_learning_diff}
\end{multline}
Note that \cref{eqt:backward_learning_diff_bound} relates the difference in the backward process approximations at the \( n \)-th iteration to the analogous difference in the forward processes also at the \( n \)-th iteration. This is in contrast to \cref{eqt:forward_learning_diff} which expresses the difference at the \( (n+1) \)-th iteration to differences at the \( n \)-th iteration. 

Defining the process 
\begin{equation}
D_{n} \coloneq \mathbb{E}(\lvert \widehat{X}_{n} - \widehat{X}_{n}^{\theta} \rvert^2),
\end{equation}
then by taking \( \mathbb{E}(\lvert \cdot \rvert^2) \) of \cref{eqt:forward_learning_diff} we obtain through Jensen's inequality,  \cref{asmp:lipschitz_continuity}, and \( \Delta t^2 \leq \Delta t \) that
\begin{align}
D_{n+1} & \preceq D_n (1 + \Delta t) + \mathbb{E}(\lvert \widehat{Y}_n - \widehat{Y}_n^{\theta} \rvert^2) \Delta t \\
\shortintertext{taking \( \mathbb{E}(\lvert \cdot \rvert^2) \) of \cref{eqt:backward_learning_diff_bound} gives}
& \preceq D_n (1 + \Delta t) + (\epsilon_{\theta}^2 + D_n) \Delta t \\
& \preceq \epsilon_{\theta}^2 \Delta t + (1 + \Delta t) D_n \\
\shortintertext{which expands recursively (with \( D_0 = 0 \)) to}
& \preceq n \epsilon_{\theta}^2 \Delta t + \Delta t \sum_{k=1}^{n} D_k \\
\shortintertext{as \( n \Delta t \leq N\Delta t = T \)}
& \preceq \epsilon_{\theta}^2 + \Delta t \sum_{k=1}^{n} D_k \\
\shortintertext{using Gr\"{o}nwall's inequality} 
& \preceq \epsilon_{\theta}^2. \label{eqt:bound_on_forward_difference_l2_strong_error}
\end{align}

We can then combine  \cref{eqt:bound_on_forward_difference_l2_strong_error} with \cref{eqt:backward_learning_diff_bound} to similarly obtain
\begin{equation}
\mathbb{E}(\lvert \widehat{Y}_{n} - \widehat{Y}_{n}^{\theta} \rvert^2) \preceq \epsilon_{\theta}^2,
\end{equation}
which is trivially combined with \cref{eqt:bound_on_forward_difference_l2_strong_error} to give the desired result, completing the proof. \qedhere
\end{proof}

\begin{corollary}
\label{thm:strong_convergence_between_em_and_nn_shrinking_tol}
If \( \epsilon_{\theta} =  O(\Delta t^{1/2}) \) then  \begin{equation}
\mathbb{E}(\lvert \widehat{X}_n - \widehat{X}_n^{\theta} \rvert^2 + \lvert \widehat{Y}_n - \widehat{Y}_n^{\theta} \rvert^2) \preceq \Delta t.
\end{equation}
\end{corollary}

\begin{remark}
Setting \( \epsilon_{\theta} =  O(\Delta t^{1/2}) \) seems an unrestrictive training criteria, that is likely satisfied by our loss function's exact minimisers already. However, we don't achieve the exact minimisers, only those we iterate to until our tolerance is achieved (or a maximum number of iterations is reached). 
\end{remark}

\begin{remark}
The temporal discretisation used in the construction of the loss function needn't match a given multilevel Monte Carlo level's temporal discretisation. Furthermore, the loss function can keep the same granularity between different levels. It is a subtle but easy mistake to conflate the discretisation used in training with that used in inference, and indeed the two may well differ. 
\end{remark}

\begin{theorem}
\label{thm:strong_convergence_between_exact_and_nn}
Under the same assumptions as \cref{thm:strong_convergence_between_em_and_nn}   \begin{equation}
\mathbb{E}(\lvert X_{n} - \widehat{X}_n^{\theta} \rvert + \lvert Y_{n} - \widehat{Y}_n^{\theta} \rvert) \preceq \max\{\epsilon_{\theta}, \Delta t^{1/2}\}.
\end{equation}
\end{theorem}

\begin{proof}
By a straightforward application of the triangle inequality we readily obtain
\begin{align}
\mathrlap{\mathbb{E}(\lvert X_{n} - \widehat{X}_n^{\theta} \rvert + \lvert Y_{n} - \widehat{Y}_n^{\theta} \rvert)} \hspace{4em}&  \notag \\
& \leq 
\begin{aligned}[t]
& \mathbb{E}(\lvert X_{n} - \widehat{X}_n \rvert + \lvert Y_{n} - \widehat{Y}_n \rvert)  \\
& {} + \mathbb{E}(\lvert \widehat{X}_n - \widehat{X}_n^{\theta} \rvert + \lvert \widehat{Y}_n - \widehat{Y}_n^{\theta} \rvert) 
\end{aligned} \\
\shortintertext{using \cref{thm:strong_error_exact_em,thm:strong_convergence_between_em_and_nn}}
& \preceq \Delta t^{1/2} + \epsilon_{\theta} 
\preceq \max \{\Delta t^{1/2}, \epsilon_{\theta} \}, 
\end{align}
which completes the proof. \qedhere 
\end{proof}

\begin{corollary}
\label{thm:strong_convergence_between_exact_and_nn_backward_process}
\begin{equation}
\mathbb{E}(\lvert Y_{n} - \widehat{Y}_n^{\theta} \rvert) \leq \max\{\epsilon_{\theta}, \Delta t^{1/2}\}.
\end{equation}
\end{corollary}

\subsection{Numerical results}
\label{sec:numerical_results}

\begin{figure*}[h!tb]
\centering
\includegraphics{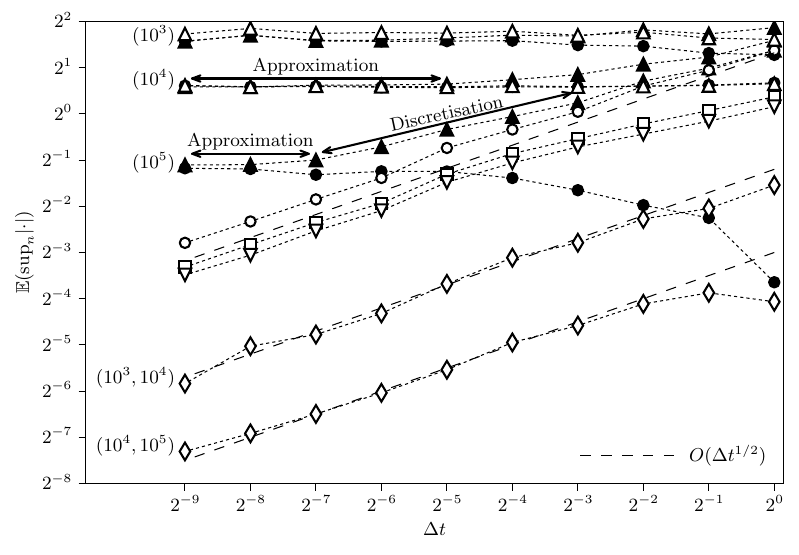}

\begin{tabular}{clcl}
%Marker & Quantity \\
%\hline
\(\CIRCLE\) & \(u(t_n, \widehat{X}_n^{\mathrm{c}}) - \hat{u}(t_n, \widehat{X}_n^{\mathrm{c},\theta};\theta)\) & \(\triangledown\) & \(\hat{u}(t_n, \widehat{X}_n^{\mathrm{f},\theta};\theta) - \hat{u}(t_n, \widehat{X}_n^{\mathrm{c},\theta};\theta)\) \\
\(\Circle\) & \(u(t_n, X_{t_n}) - u(t_n, \widehat{X}_n^{\mathrm{c}})\) & \(\blacktriangle\) & \(u(t_n, X_{t_n}) - \hat{u}(t_n, \widehat{X}_n^{\mathrm{c},\theta};\theta)\) \\
\(\Square\) & \(u(t_n, \widehat{X}_n^{\mathrm{f}}) - u(t_n, \widehat{X}_n^{\mathrm{c}})\) &
\(\vartriangle\) & \(\hat{u}(t_n, \widehat{X}_n^{\mathrm{f},\theta'};\theta') - \hat{u}(t_n, \widehat{X}_n^{\mathrm{f},\theta};\theta)\) \\
& &
(\(\vartriangle\)) & \(\hat{u}(t_n, \widehat{X}_n^{\mathrm{f},\theta'};\theta') - \hat{u}(t_n, \widehat{X}_n^{\mathrm{c},\theta};\theta)\) \\
\(\amsDiamond\) & \multicolumn{3}{l}{\(\hat{u}(t_n, \widehat{X}_n^{\mathrm{f},\theta'};\theta') - \hat{u}(t_n, \widehat{X}_n^{\mathrm{c},\theta'};\theta') - \hat{u}(t_n, \widehat{X}_n^{\mathrm{f},\theta};\theta) + \hat{u}(t_n, \widehat{X}_n^{\mathrm{c},\theta};\theta)\)}
\end{tabular}
\caption{The strong error using the \(L^1\)-norm from \cref{eqt:strong_error_lp} for various two-way and four-way differences of significance to multilevel Monte Carlo settings. Reference lines for a strong convergence order of \(\tfrac{1}{2}\) are shown. Terms using a trained model parametrised by some \(\theta\) have the number of iterations used to train the model shown in parentheses (e.g.\ annotating the \(\blacktriangle\)~markers). Similarly, terms using the same model at two different stages in training with differing numbers of iterations performed are denoted by \(\theta\) and \(\theta'\), and both the respective iterations are parenthesised (e.g.\ annotating the \(\amsDiamond\)~markers). For some example models, we indicate regions where the relevant error is dominated either by the discretisation error, or the approximation error for an imperfectly trained model. A legend marker which is parenthesised indicates that the results map directly onto those of the unparenthesised marker, but are omitted for clarity.}
\label{fig:multilevel_variances}

\vspace{1\baselineskip}
\begin{tabular}{llc}
Result & Summary & Marker \\
\hline
\firstrowspacing
\Cref{thm:strong_error_decoupled_forward,thm:strong_error_exact_em} & \( \mathbb{E}(\lvert Y_n - \widehat{Y}_n \rvert) \preceq \Delta t^{1/2} \) & \(\Circle\) \\
\Cref{thm:strong_convergence_between_em_and_nn} & \( \mathbb{E}(\lvert \widehat{Y}_n - \widehat{Y}_n^{\theta} \rvert) \preceq \epsilon_{\theta} \) & \(\CIRCLE\) \\
\Cref{thm:strong_convergence_between_exact_and_nn_backward_process} & \( \mathbb{E}(\lvert Y_n - \widehat{Y}_n^{\theta} \rvert) \preceq \max \{\epsilon_{\theta}, \Delta t^{1/2} \} \) & \(\blacktriangle\)
\end{tabular}
\captionof{table}{Summary of the expected bounding behaviours arising from our numerical analysis, and the corresponding markers in \cref{fig:multilevel_variances}.}
\label{tab:expected_multilevel_variances}
\end{figure*}

We can now assess whether the strong error bounds derived thus far appear in practice, (the most significant of which are \cref{thm:strong_error_exact_em,thm:strong_convergence_between_em_and_nn,thm:strong_convergence_between_exact_and_nn_backward_process}). Our corresponding numerical results are shown in \cref{fig:multilevel_variances}, and in \cref{tab:expected_multilevel_variances} we surmise our main analytic results and how these relate to the corresponding data shown in \cref{fig:multilevel_variances}.

Considering first \cref{thm:strong_error_decoupled_forward}, and its more general \( L^2 \)-norm extension for coupled processes,   \cref{thm:strong_error_exact_em}. In \cref{fig:multilevel_variances} this corresponds to the data for \( u(t_n, X_{t_n}) - u(t_n, \widehat{X}_{n}) \) (\(\Circle\)~markers) and predicts a strong decay rate of \( \Delta t^{1/2} \), which is exactly what we see in \cref{fig:multilevel_variances}. Overall, this result is not surprising, as it is in keeping with the usual equivalent results obtained for the forward processes, and also the similar results from \citet{negyesi2024generalizedconvergencedeepbsde}. 

From \cref{thm:strong_convergence_between_em_and_nn} we would anticipate the data for \( u(t_n, \widehat{X}_{n}) - \hat{u}(t_n, \widehat{X}_{n}^{\theta};\theta) \) (\(\CIRCLE\)~markers) to exhibit a constant value (of order \( \epsilon_{\theta} \)). For the models trained with \SI{e3}{} and \SI{e4}{} iterations, this is exactly what we observe. Interestingly, while we do observe the same trend for the same model trained with \SI{e5}{} iterations, it does show spuriously better convergence for the coarsest few time steps. A heuristic explanation for this can be garnered by examining the results from \citet{raissi2018forwardbackward}, whereby their model does not show a particularly uniform accuracy, but instead shows a greater relative accuracy at the starting and terminal points.
Consequently, coarser grids will by construction have the initial and terminal errors predominantly constitute the set of errors being evaluated in the supremum, and thus leads to the considerable reduction in the strong error. 

Lastly, from \cref{thm:strong_convergence_between_exact_and_nn_backward_process} we would expect \( u(t_n, X_{t_n}) - \hat{u}(t_n, \widehat{X}_{n}^{\theta};\theta) \) (\(\blacktriangle\)~markers) to demonstrate a strong decay rate which transitions from a \( \Delta t^{1/2} \) rate for coarse time steps to a constant \( \epsilon_{\theta} \) value for sufficiently fine time steps, which is exactly what we observe in \cref{fig:multilevel_variances}. The transition from discretisation dominated into approximation dominated is highlighted in the figure, and the onset is later for better trained models, as expected. 

In our proof of the results supporting  \cref{thm:strong_convergence_between_exact_and_nn_backward_process} we made extensive use of
Jensen's inequality and Gr\"{o}nwall's inequality. Neither of these should be expected to produce particularly tight bounds, and thus we do not expect that had we kept closer account of all the numerous coefficients (absorbed into our ``\( \preceq \)'' notation from \cref{def:relational_o_notation}) that  \cref{thm:strong_convergence_between_exact_and_nn_backward_process} would be a tight bound. So while \cref{thm:strong_convergence_between_exact_and_nn_backward_process} is theoretically useful, it is not expected to be a practical tool to determine when we transition from discretisation to approximation dominated regimes. To emphasise this point, we quote \citet[p.\,7]{iserles2009first} who describes analogous bounds arising from applications of Gr\"{o}nwall's inequality for analysing ordinary differential equations: ``\textit{At first sight, it might appear that there is more to the last theorem than meets the eye---not just a proof of convergence but also an upper bound on the error. In principle this is perfectly true [\ldots] The problem with this bound is that, unfortunately, in an overwhelming majority of practical cases it is too large by many orders of magnitude [\ldots] The moral of our discussion is simple. The bound from [the proof] must not be used in practical estimations of numerical error!}''. This shortcoming is equally true for our bounds. Finding such a practical estimator remains a topic for further research.

While our focus has been on bounding the leading order two-way differences, we can recall the four-way difference from \cref{eqt:nested_mlmc_nn_and_temporal}, for which we have shown the corresponding results also in \cref{fig:multilevel_variances} (\(\amsDiamond\)~markers). Indeed we can see that the variance is considerably lower than the other terms, and shows a strong convergence order of \( \tfrac{1}{2} \). Both these behaviours are encouraging for multilevel Monte Carlo applications. Nonetheless, it is the intermediary terms in \cref{eqt:nested_mlmc_nn_and_temporal} and the correction term in \cref{eqt:mlmc_nn_and_temporal} which dominate, and which would benefit from further analytic attention. 
\section{Future research}
\label{sec:future_research}

Over the course of this report we have taken the setup popularised by \citet{e2017deep_learning,e2018solving,e2021algorithms} and \citet{raissi2018forwardbackward}. However, there have arisen topics which merit consideration for further investigation. Some we have mentioned thus far when they first arose, and others we newly discuss here. 

\subsection{The loss function's bias and variance}

We discussed in detail the bias and variance properties of the loss function. Through our discussions we indicated the optional inclusion of the additional term \cref{eqt:loss_extra_term}. However, the benefits from including (or excluding) this term appear to have received insufficient attention. Its use seems to have proliferated into various implementations without sufficient evidence in regards to its impact. Also, whether it should be equally weighted or not has also received insufficient attention. 
Similarly, we suggested the further inclusion of a higher order term resulting in 
\cref{eqt:loss_higher_order}. The impact of this and the [in]significance of the costs from computing the Hessian should be investigated. Lastly, in \cref{remark:new_brownian_motion_samples_between_iterations}, and again when discussing \cref{algo:euler_maruyama_with_nn} on page~\pageref{paragraph:loss_function_resampling}, we raised the issue of whether the Brownian paths should be resampled between different training iterations, or kept constant. This finer detail, and the impact this has on convergence in training and also on the bias in the loss function would benefit from further research. 

The loss function \cref{eqt:loss} was strongly advocated for and praised by \citet{raissi2018forwardbackward}, discussing the advantages over using the more conventional and wider spread loss functions which only measure convergence at the terminal condition for the backward process. However, when we compare our bounds to the equivalent counterparts from
e.g.\ \citet{negyesi2024generalizedconvergencedeepbsde} (which we restated in \cref{thm:negyesi2024generalizedconvergencedeepbsde}), we see both methods have obtained the same order of strong convergence. This puts into question the supposed benefits of the loss function \cref{eqt:loss}, at least from a theoretical viewpoint. Consequently, a more rigorous benchmark comparison between the different formulations of the loss function certainly seems worthwhile. 

Lastly, recognising the two differing complexities of the candidate loss functions, they present another form of possible multilevel decomposition. A crude estimator formed using a cheap and quick loss function, and a fine estimator using a more sophisticated but expensive loss function. They could differ in which terms are included or excluded (e.g.\ such as whether \cref{eqt:loss_extra_term} is added to \cref{eqt:loss}), or by their batch sizes and discretisation levels.

\subsection{Variance reduction techniques}

When having discussions of reducing the variance of quantities involved in Monte Carlo simulations, there is a wealth of literature on variance reduction techniques for estimators, as discussed by \citet[\S\,4]{glasserman2003monte} and \citet[\S\,V]{asmussen2007stochastic}. The Brownian motion paths used to define the loss function we asserted were independently and identically distributed sample paths. However, to reduce the variance introduced by our Brownian motion sample paths, we might consider similarly coupling these. 

The most immediately applicable variance reduction technique would be to use \textit{antithetic variates}. For a comprehensive detailing, we recommend the reader to \citet[\S\,4.2]{glasserman2003monte}. In summary, for Wiener increments \( \Delta W_n \coloneq \sqrt{\Delta t} Z_n \coloneq \sqrt{\Delta t} \Phi^{-1}(U_n)\), where \( Z_n \) is a standard Gaussian random variable, which can be formed by the inverse transform method \citep{glasserman2003monte} from a standard uniform random variable \( U_n \) using the standard Gaussian inverse cumulative distribution function \( \Phi^{-1} \), we also generate a second path produced with the underlying uniform random variables \( 1 - U_n \), producing a reflected Brownian motion path. 

Another would be to use the \textit{antithetic twin} paths proposed by \citet{giles2014antithetic}, whereby the fine increments \( \Delta W_n^{\mathrm{f}} \) and \(\Delta W_{n + 1/2}^{\mathrm{f}} \) in the Euler-Maruyama approximation for the fine path \cref{eqt:euler_maruyama_fine} are swapped in the order they are used. 

Considering these methods, they appear readily applicable for constructing a modified version of the loss function in \cref{eqt:loss} which could utilise either approach (or both). Notably, the approach using antithetic twin paths has been used to circumvent the need for simulating L\'{e}vy areas, and thus may be appropriate should researchers in the future wish to utilise higher order numerical schemes beyond the Euler-Maruyama scheme. 

\subsection{Better interpolation points}

A remark we find interesting is that the loss function in \cref{eqt:loss} measures the loss at equally spaced time points for times in the domain \( [0, 1] \). For the analysis of the Euler-Maruyama scheme, it is extremely convenient to assume the time steps are of a uniform size. In general though, we usually only require that \( \lim_{N\to \infty} \sup_{n < N} \Delta t_n = 0 \), and we can allow for unequal time steps. The loss function we notice looks to learn an approximation which minimises the error at these equally spaced points, and can be interpreted as trying to learn an approximation that can interpolate between these points and has zero error at these interpolation points. From the viewpoint of a numerical analyst, choosing the points from which to form an interpolation is a well studied topic, and the answer in practice is well known to never use equally spaced intervals. Instead, points should be sampled more densely at the ends of the interval of interest, and the example \textit{par excellence} of good interpolation points are Chebyshev points. For excellent resources on this subject, we recommend the reader to \citet{trefethen2019approximation} and \citet{powell1981approximation}. Whether learning a loss function evaluated at such interpolation points would yield a superior results we believe would be an interesting topic for further investigation. 

The results from \citet{raissi2018forwardbackward} show the relative error being quite uniform over the entire time domain, and in fact noticeably better at the starting and terminal times. This is in contrast to regular interpolation errors on equispaced intervals which grow wildly near the edges \citep{trefethen2019approximation,powell1981approximation}. One may then have the impression that the suggested research into non-equidistant interpolation points would be moot. Such a conclusion though would be misleading, as the errors shown  by \citet{raissi2018forwardbackward} are measured at the interpolation points, not between them. Consequently, we can't draw conclusions about errors between the interpolation points (and in regular numerical analysis is it between the interpolation points where the errors can grow horribly, e.g.\ Runge phenomenon). One possible benefit from such research might be that we can recover bounds on the uniformity of the convergence of the neural network approximation, whereas we have thus far had to make additional assumptions (i.e.\ \cref{asmp:uniform_convergence}). We recall the uniformity results from \citet{perotti2024pricing} (which we restated in \cref{thm:perotti2024pricing}), which would be very useful to parallel in our setup.

\section{Conclusions}
\label{sec:conclusions}

Having first briefly motivated continuously evolving stochastic systems in \cref{sec:introduction}, we provided in \cref{sec:mathematical_preliminaries} a comprehensive overview of all the mathematical preliminaries and requisites that encompass the coupling of multilevel Monte Carlo methods with neural networks for simulating forward-backward stochastic differential equations. Not only did this lay the groundwork within which we could subsequently present our research, it further allowed us to highlight the differing setups between isolated online and offline training. Additionally, we were able to posit appropriate multilevel Monte Carlo decompositions which incorporated neural networks. 

With the broad body of background material presented, we focused on the most directly neighbouring works of research closest to our own \citep{e2017deep_learning,e2018solving,e2021algorithms,raissi2018forwardbackward,guler2019robust,naarayan2024thesis} in \cref{sec:previous_research}. The works of \citet{raissi2018forwardbackward} and \citet{guler2019robust} were detailed because of the specific loss function advocated by \citet{raissi2018forwardbackward}, in contrast to more usual formulations (e.g.\ \citep{e2017deep_learning,e2018solving,e2021algorithms}). We inspected the experimental framework proposed by \citet{guler2019robust}, presenting \cref{algo:mlmc_inspired_training} as a formalisation of their ``multilevel Monte Carlo inspired'' training framework. We discussed that the original ``multilevel Monte Carlo'' framework proposed by \citet{guler2019robust} was a misnomer, and instead offered alternatives which respected the telescoping summation central to multilevel Monte Carlo, and detailed the multilevel coupling mechanisms we built into \cref{algo:mlmc_inspired_training}.

The setup by \citet{raissi2018forwardbackward} convincingly extols the benefits of a path wise loss function construction (cf.\ \cref{eqt:loss}) through impressive empirical results. However, while the empirical results are striking, there is no supporting analysis to accompany the proposed setup. This is both displeasing from a theoretically aesthetic perspective, but also limiting for practitioners wishing to adopt such frameworks in multilevel Monte Carlo setups. To address this gap, in \cref{sec:numerical_analysis} we scrutinised the loss function proposed by \citet{raissi2018forwardbackward}, giving a heuristic analysis which quantified the bias and variance of the proposed loss function. From this analysis, we were able to propose alternative loss function formulations which we hope should exhibit reduced variance (by a factor of \( \Delta t^{1/2} \)). Thereafter, we produced novel strong error bounds for numerical approximations of coupled forward-backward stochastic differential equations utilising neural network approximations using \cref{algo:euler_maruyama_with_nn}, with supporting experimental results. These analytic bounds closely resembled the usual counterparts in the existing literature, namely the classic strong convergence results presented by \citet{kloeden1999numerical}, and also very recent analogous results by \citet{negyesi2024generalizedconvergencedeepbsde}, who use differing problem setups which also utilise neural networks.
While the bounds we produced are theoretically useful, they are not tight. This means more work is required for fully prescriptive tools which can determine whether temporal discretisations need to be refined, or neural networks trained further or extended to greater numbers of layers and neurons.\footnote{Semi-prescriptive bounds on the neural network's shape are given by \citet{perotti2024pricing}.} 

Having covered the aforementioned topics, we highlighted various avenues for further research in \cref{sec:future_research}, where all the topics centred around the loss function. Concerning issues touched upon during the course of our research was whether uniform convergence bounds (such as similar bounds by \citet{perotti2024pricing}) could be achieved, and if batches of Brownian paths should be resampled between training iterations. In close proximity to our research was whether antithetic techniques could be put to good use in our problem setup. As a last topic closer to classical numerical analysis was whether non-equidistant interpolation points could be used (e.g.\ Chebyshev points), and if these might give rise to uniform convergence bounds that we have had to otherwise assume.

%\newpage
\begin{footnotesize} % footnotesize, small, normalsize
\bibliography{references}
\end{footnotesize}

\end{document}